\documentclass[12pt, reqno]{amsart}

\allowdisplaybreaks

\usepackage[all]{xy}
\usepackage{amssymb} 

\newcommand{\morespace}{}

\newcommand{\wt}{\widetilde}
\newcommand{\wh}{\widehat}
\newcommand{\ov}{\overline}

\newcommand{\mybeginexample}[1]{\vspace{0.1 true cm}\emph{Example #1.}}
\newcommand{\myendexample}{\vspace{0.1 true cm}}

\newcommand{\ii}{\mathrm{i}}
\newcommand{\supp}{\mathrm{supp}}
\newcommand{\natlog}{\mathrm{ln}}

\newtheorem{thm}{Theorem}
\newtheorem{prop}{Proposition}
\newtheorem{defn}{Definition}

\begin{document}
\morespace

\title[$q$-Legendre transformation]{
$q$-Legendre transformation: partition\\ 
functions and quantization\\ 
of the Boltzmann constant
}

\author[A. E. Ruuge]{Artur E. Ruuge}
\author[F. van Oystaeyen]{Freddy van Oystaeyen}

\address{
Department of Mathematics and Computer Science, 
University of Antwerp, 
Middelheim Campus Building G, 
Middelheimlaan 1, B-2020, 
Antwerp, Belgium
}
\email{artur.ruuge@ua.ac.be, fred.vanoystaeyen@ua.ac.be}

\begin{abstract}
\morespace
In this paper we construct a $q$-analogue of the Legendre transformation, where 
$q$ is a matrix of formal variables defining the phase space braidings between the 
coordinates and momenta (the extensive and intensive thermodynamic observables). 
Our approach is based on an analogy between the semiclassical 
wave functions in quantum mechanics and the quasithermodynamic partition 
functions in statistical physics. 
The basic idea is to go 
from the $q$-Hamilton-Jacobi equation in mechanics  
to the $q$-Legendre transformation in thermodynamics. 
It is shown, that this requires a non-commutative analogue of the 
Planck-Boltzmann constants ($\hbar$ and $k_{B}$) 
to be introduced back into the classical formulae. 
Being applied to statistical physics, 
this naturally leads to an idea to go further and to replace the Boltzmann constant 
with an infinite collection of generators of the so-called epoch\'e (bracketing) algebra. 
The latter is an infinite dimensional noncommutative algebra 
recently introduced in our previous work, which 
can be perceived as an infinite sequence of 
``deformations of deformations'' of the Weyl algebra. 
The generators mentioned are naturally indexed by planar binary 
leaf-labelled trees in such a way, that the trees with a single leaf correspond 
to the observables of the limiting thermodynamic system.   
\end{abstract}

\maketitle

\section{Introduction}

The Legendre transformation 
plays a rather fundamental role in mathematical physics. 
One can immediately think of two examples:  
\emph{phenomenological} thermodynamics and 
\emph{classical} mechanics. 
The free energy $F$ and the internal energy $E$ of a thermodynamic system are 
related to each other in the same way as  
the Lagrangian $L$ and the Hamiltonian $H$ of a mechanical system. 
At this level, the Legendre transformation is just a 
convenient concept that connects different pictures 
of description of a physical system.  

Much more important is that the Legendre transformation is associated to 
the transition from quantum statistical physics  to the classical limit. 
It emerges both in the semiclassical approximation of 
\emph{quantum} mechanics, 
as well as in the quasithermodynamic approximation of 
\emph{statistical} physics. 
In the first case, one needs to 
consider the (additive) asympotics of the wave 
functions of the coordinate and momentum representations 
corresponding to $\hbar \to 0$, 
and in the second case one needs to 
consider the (multiplicative) asymptotics of the 
partition functions of the canonical and 
microcanonical ensembles corresponding to $k_{B} \to 0$. 
There is an analogy between these two limits  
corresponding to a similarity between 
fast oscillating wave functions and rapidly decaying partition functions.

The original motivation of this paper was to study a $q$-deformation of the 
Legendre transformation (where $q$ is a matrix of formal variables) 
from the perspective of these limit transitions. 
It turns out, that the consequences of these investigations 
seem to go much deeper than one would expect and indicate a necessity to replace 
the fundamental constants $\hbar$ and $k_{B}$ with \emph{infinite} collections of 
quantities $\lbrace \hbar_{\Gamma} \rbrace_{\Gamma}$ and $\lbrace (k_{B})_{\Gamma} \rbrace_{\Gamma}$, 
respectively. 
These collections correspond to the generators of the 
\emph{epoch\'e} algebra introduced in \cite{RuugeVanOystaeyen}, and 
are structured in a slightly more complicated way than infinite dimensional 
matrices (the Heisenberg's quantization). The index $\Gamma$ varies over 
a set of finite leaf-labelled planar binary trees of different sizes, 
while the labelling set is just the set of symbols 
corresponding to the degrees of freedom of the 
limiting physical system. 
We start with describing the four examples mentioned in a little more detail. 

\mybeginexample{1}
If $F = F (T, V)$ is the free energy of a 2-dimensional thermodynamic system at 
absolute temperature $T$ and macroscopic volume $V$, then the corresponding entropy 
$S$ is given by $S = - (\partial F/ \partial T)_{V}$, and the 
corresponding pressure is given by $p = - (\partial F/ \partial V)_{T}$. 
Taking the Legendre transformation of $F (T, V)$ with respect to $T$ yields the 
internal energy $E$ of the system in terms of $S$ and $V$, $E = E (S, V)$, and the 
Legendre transformation with respect to $V$ defines the Gibbs potential $G (T, p)$. 
The restriction of the 1-form 
\begin{equation*} 
\alpha := T d S - p d V
\end{equation*}
to the submanifold of the equilibrium states $\Lambda \subset \mathbb{R}^{4} (S, V; T, p)$
is exact: $\alpha|_{\Lambda} = d E|_{\Lambda}$, since 
$(T d S)|_{\Lambda} = \delta Q$ (the infinitesimal amount of heat 
absorbed by the system) and 
$(p d V)|_{\Lambda} = \delta A$ (the infinitesimal work produced by the system) 
add up according to the law of the conservation of energy 
to $d E|_{\Lambda} = \delta Q - \delta A$. 
\myendexample

\mybeginexample{2}
This basic example of the Legendre transformation stems from \emph{classical} mechanics. 
If $L = L (x, v)$ is the Lagrangian of a system with $d$ degrees of freedom 
described by coordinates $x \in \mathbb{R}^{d}$ and the 
associated velocities $v \in \mathbb{R}^d$, then one may consider 
the momenta 
$p_i = \partial L (x, v)/ \partial v_i$, $i = 1, 2, \dots, d$. 
Assuming that this implicitly defines $v = \wt{v} (x, p)$, $p = (p_1, p_2, \dots, p_d)$, 
one can switch to the Hamiltonian formalism via the Legendre transformation 
\begin{equation*} 
H (x, p) = (p v - L (x, v))|_{v = \wt{v} (x, p)},  
\end{equation*}
where $p v := \sum_{i = 1}^{d} p_i v_i$. 
Introducing $f_i := \partial L (x, v)/ \partial x_{i}$, $i = 1, 2, \dots, d$, 
one may interpret the states of the system as a $2 d$-dimensional (Lagrangian) 
submanifold $\Lambda^{\#} \subset \mathbb{R}^{4 d} (p, x; v, f)$, $f = (f_1, f_2, \dots, f_d)$, 
such that the restriction of the 1-form 
\begin{equation*} 
\varepsilon := \sum_{i = 1}^{d} (v_i d p_i - f_i d x_i) 
\end{equation*}
to $\Lambda^{\#}$ is exact, 
the conservation of energy being nothing else but 
$\varepsilon|_{\Lambda^{\#}} = d H|_{\Lambda^{\#}}$. 
\myendexample

The Legendre transformation 
becomes conceptually more important if one takes 
a step from classical mechanics to \emph{quantum} mechanics and from 
the phenomenological thermodynamics to the \emph{statistical} physics of equilibrium states. 
It is quite remarkable that 
these two steps are associated with a pair of fundamental physical constants 
\begin{equation*} 
\begin{aligned}
\hbar &= 1.05 \times 10^{-27} \mathit{erg} \, s, \\
k_{B} &= 1.38 \times 10^{-16} \mathit{erg} \, K^{-1}, 
\end{aligned}
\end{equation*}
the Planck's constant $\hbar$ and the Boltzmann constant $k_{B}$. 
Naively, $\hbar$ corresponds to the quantization of 
the spectra of ``physical properties'', and $k_{B}$ 
corresponds to the quantization of the ``chemical substance'' itself. 
Let us consider the transition from quantum mechanics to classical mechanics. 
This is equivalent to pretending that $\hbar \to 0$ 
(in fact, this will be a dimensionless combination of $\hbar$ and some scaling parameters 
defining the natural units of measurement in the experimental set-up). 
To avoid a confusion in what follows, one may wish to denote the physical values of 
$\hbar$ and $k_{B}$ mentioned above as $\hbar^{\mathit{phys}}$ and $k_{B}^{\mathit{phys}}$.  

\mybeginexample{3}
Suppose the Schr\"odinger equation for a given physical system is of the shape 
\begin{equation*} 
\ii \hbar \frac{\partial \psi_{\hbar}^{t} (x)}{\partial t} = 
H \Big( \overset{2}{x}, - \ii \hbar \frac{\overset{1}{\partial}}{\partial x} \Big) \psi_{\hbar}^{t} (x),  
\end{equation*}
where $x \in \mathbb{R}^d$ are the classical coordinates, 
$t \in \mathbb{R}$ is classical time, 
$\psi_{\hbar}^{t} (x) \in L^2 (\mathbb{R}^d)$ is the wave-function, 
$H (x, p)$ is (for simplicity) a polynomial in $x$ and $p \in \mathbb{R}^d$ 
(the classical momenta), and the indices atop denote the order of action. 
If we switch the representation of the algebra of observables to an isomorphic one 
by taking the $\hbar$-Fourier transform, 
$\wt{\psi}_{\hbar}^{t} (p) = (2 \pi \hbar)^{-d /2} 
\int dx \, \exp (- \ii p x /\hbar) \psi_{\hbar}^{t} (x)$, then 
the shape of the Schr\"odinger equation is just as good: 
\begin{equation*} 
\ii \hbar \frac{\partial \wt{\psi}_{\hbar}^{t} (p)}{\partial t} = 
H \Big( \ii \hbar \frac{\overset{2}{\partial}}{\partial p}, \overset{1}{p} \Big) \wt{\psi}_{\hbar}^{t} (p). 
\end{equation*}
Suppose $\psi_{\hbar}^{t} (x)$ is described for some $t \in [0, T]$ by an additive asymptotics 
$\psi_{\hbar}^{t} (x) = \varphi^{t} (x) \exp (i S^t (x)/ \hbar) / \sqrt{J^t (x)} + O (\hbar)$, 
where $\hbar \to 0$, $S^t (x)$ is a real-valued smooth function in $x$ and $t$, 
$\varphi^{t} \in C_{0}^{\infty} (\mathbb{R}^d)$, for each $t \in [0, T]$, and 
$J^t (x) = |\det S_{x x}^{t} (x)| \not = 0$, for $t \in [0, T]$ and 
$x \in \supp \varphi^{t}$, where 
$S_{x x}^{t} (x) := \| \partial^2 S^t (x) / \partial x_i \partial x_j\|_{i, j = 1}^{d}$ is the Hess matrix of 
$S^t (x)$. Then $\wt{\psi}_{\hbar}^{t} (p) = \wt{\varphi}^{t} (p) 
\exp (i \wt{S}^t (p)/ \hbar) / \sqrt{\wt{J}^t (p)} + O (\hbar)$, 
where the functions $\wt{S}^t (p)$, $\wt{J}^t (p)$, and $\wt{\varphi}^{t} (p)$ 
are smooth and can be computed using the 
stationary phase method, 
\begin{equation*} 
\frac{1}{(2 \pi \hbar)^{d/ 2}}
\int dx\, \exp \Big( - \frac{\ii p x}{\hbar} + \frac{\ii S^t (x)}{\hbar} \Big) 
\frac{\varphi^t (x)}{\sqrt{J^t (x)}} = 
\exp \Big( \frac{\ii \wt{S}^t (p)}{\hbar} \Big)
\frac{\wt{\varphi}^{t} (p)}{\sqrt{\wt{J}^t (p)}} + O (\hbar), 
\end{equation*}
where $O (\hbar)$ are the terms of the formal asymptotic expansion in $\hbar \to 0$ of order $\geqslant 1$, 
$\wt{J}^t (p) = |\det \wt{S}_{p p}^{t}|$, 
$\wt{S}_{p p}^{t} := \| \partial^2 \wt{S}^t (p)/ \partial p_{i} \partial p_{j} \|_{i, j = 1}^{d}$, 
the radical sign denotes the principal square root. 
The function $\wt{S}^t (p)$ is just the Legendre transformation of $S^t (x)$, 
\begin{equation*} 
\wt{S}^t (p) = (- p x + S^t (x))|_{x = \wt{x}^t (p)}, 
\end{equation*} 
where $x = \wt{x}^t (p)$ is the solution of $p = \frac{\partial S^t (x)}{\partial x}$ 
with respect to $x = (x_1, x_2, \dots, x_d)$ if $p = (p_1, p_2, \dots, p_d)$ is perceived as a parameter. 
To link the functions $\varphi^t (x)$ and $\wt{\varphi}^t (p)$, it is convenient 
to consider the graph 
$\Lambda^t := \lbrace (x, p) | p = \partial S^t (x)/ \partial x \rbrace 
\subset \mathbb{R}_{x}^{d} \times \mathbb{R}_{p}^{d}$, 
or, what is the same, 
$\Lambda^t = \lbrace (x, p) | x = - \partial \wt{S}^t (p)/ \partial p \rbrace 
\subset \mathbb{R}_{x}^{d} \times \mathbb{R}_{p}^{d}$, 
and then lift the both functions to $\Lambda^t$ with respect to the 
canonical projections $\pi_x: \mathbb{R}_{x}^{d} \times \mathbb{R}_{p}^{d} \to \mathbb{R}_{x}^{d}$ and  
$\pi_p : \mathbb{R}_{x}^{d} \times \mathbb{R}_{p}^{d} \to \mathbb{R}_{p}^{d}$, respectively. 
Denoting the result corresponding to $\varphi^t$ as $\chi^t \in C_{0}^{\infty} (\Lambda^t)$, 
and the result corresponding to $\wt{\varphi}^t$ as $\wt{\chi}^t \in C^{\infty} (\Lambda^t)$, 
one can check using the explicit formulae of the stationary phase method, that 
$\chi^t$ and $\wt{\chi}^t$ coincide on  
$U^t := (\pi_x|_{\Lambda^t})^{-1} (\supp \varphi^t)$ 
up to a phase factor. More precisely, 
\begin{equation*} 
\chi^t|_{U^t} = \exp \Big( \frac{\ii \pi}{4} M \Big) \wt{\chi}^t|_{U^t}, 
\end{equation*}
where $M \in \mathbb{Z}$. 
Computing this integer is actually quite important, 
since it leads to the discovery of the \emph{Maslov index}. 
Up to this point, we have the third example of the Legendre transformation 
stemming from the semiclassical limit $\hbar \to 0$ of quantum mechanics.  
Informally, the Legendre transformation shows up as a fast oscillating limit 
of the Fourier transform. 
It is interesting to mention that it can be 
perceived as an idempotent analogue of the Fourier transform
in the framework of ``idempotent'' functional analysis and ``tropical'' algebraic geometry 
\cite{LitvinovMaslovShpiz} 
if one replaces the usual integrals $\int$ with the idempotent integrals $\int_{\oplus}$. 
\myendexample

\mybeginexample{4}
Finally, let us consider the fourth example of the Legendre transformation. 
In analogy with the semiclassical approximation of quantum mechanics where one deals with $\hbar \to 0$, 
the transition from the statistical physics of equilibrium states can be formally understood 
as a limit $k_{B} \to 0$. 
Indeed, consider a system placed in a thermostat. Its equilibrium state  
is described by the partition function 
\begin{equation*} 
Z_{E_{*}} (\beta) := \sum_{m} \exp \Big( - \frac{1}{k_B} \beta E_{m} \Big), 
\end{equation*}
where $m$ is the index of a microscopic state, $E_{m}$ is the corresponding energy, 
$E_{*} := (E_{m})_{m}$, 
the sum is taken over all microscopic states, and 
$\beta = 1/ T$ is the inverse temperature of the thermostat in the absolute thermodynamic scale. 
In the thermodynamic limit (i.e. the macroscopic volume $V \to \infty$, 
while the value of the density of mass $\rho > 0$ is a fixed positive constant),  
for every fixed $\beta = 1/ T$, we have 
$\natlog Z_{E_{*}} (\beta) = O (N)$, where 
$N \to \infty$ is the number of corpuscles in the system.  
If we choose a 
unit of measurement $\varepsilon_{0}$ of the specific energy,   
then there is a small parameter 
$k' := k_{B} T/ (N \varepsilon_{0})$ 
(see \cite{Maslov1, MaslovMishchenko, MaslovNazaikinskii}). 
It can be expressed as 
$k' = k_{B} T m_{0}/ (V \rho \varepsilon_{0})$, 
where $m_{0}$ is the corpuscular mass. 
Since the parameters $m_{0}$, $\varepsilon_{0}$, $T$, and $\rho$ 
are fixed under the limit transition, 
essentially one deals with $k_{B}/ V$.  
Therefore, $k_{B} \to 0$, $V$ -- fixed, is ``the same thing'' 
as $k_{B}$ -- fixed, $V \to \infty$.  
The thermodynamic limit $k' \to 0$ is formally equivalent $k_{B} \to 0$. 
Now, assuming the spectrum of energies is bounded from below and satisfies 
$\min ( \lbrace E_{m} \rbrace_{m} ) \to 0$ as $k_{B} \to 0$, 
we have:  
\begin{equation*} 
Z_{E_{*}} (\beta) = \int_{0}^{\infty} 
\exp \Big( - \frac{1}{k_{B}} \beta E \Big) 
W_{k_{B}} (E, V) d E \cdot \big\lbrace 1 + o (k_{B}) \big\rbrace, 
\end{equation*}
for some density $W_{k_{B}} (E, V)$
of the measure of integration over $d E$, 
the parameter $V$ is the macroscopic volume of the system. 
In the equilibrium statistical physics, 
the quantity $F_{E_{*}} (T) := - k_{B} \beta^{-1} \natlog Z_{E_{*}} (\beta)$, $\beta = 1/ T$,  
(termed the \emph{free energy} of the canonical Gibbs ensemble at temperature $T = \beta^{-1}$)
has a non-vanishing limit as $k_{B} \to 0$, 
$F_{E_{*}} (T) = F (T, V) \lbrace 1 + o (k_{B})\rbrace$. 
To ensure this, one takes $W_{k_{B}} (E, V) = (A (E, V) / \sqrt{2 \pi k_{B}} ) 
\exp(S (E, V)/ k_{B}) \lbrace 1 + o (k_{B}) \rbrace$, for some $A (E, V)$ and $S (E, V)$. 
Then, the saddlepoint method yields: 
\begin{multline} 
\label{eq:saddle}
\frac{1}{\sqrt{2 \pi k_{B}}}
\int_{0}^{\infty} d E \, 
\exp \Big( - \frac{1}{k_{B}} \beta E + \frac{1}{k_{B}} S (E, V) \Big) A (E, V) = \\ = 
\frac{A (E, V)}{\sqrt{- \partial^2 S (E, V)/ \partial E^2}} \bigg|_{E = \wt{E} (\beta, V)}
\exp \Big( - \frac{1}{k_{B}} \Psi (\beta, V) \Big) \cdot 
\lbrace 1 + O (k_{B}) \rbrace, 
\end{multline}
where 
$E = \wt{E} (\beta, V)$ is the solution of equation $\beta = \partial S (E, V)/ \partial E$ 
with respect to $E$ (perceiving $\beta$ and $V$ as parameters), and 
the function $\Psi (\beta, V)$ is the Legendre transformation of 
$S (E, V)$ in the variable $E$, 
\begin{equation*} 
\Psi (\beta, V) = (\beta E - S (E, V))|_{E = \wt{E} (\beta, V)}. 
\end{equation*}
Since the leading term in the right-hand side of 
the previous equation 
(\ref{eq:saddle}) 
must be just 
$\exp (- \beta F (\beta^{-1}, V)/ k_{B})$, one concludes, that 
$\Psi (\beta, V) = \beta F (\beta^{-1}, V)$, and $A (E, V) = - \partial^{2} S (E, V)/ \partial E^2$. 
One has the fourth example of the Legendre transformation. 
\myendexample

Recall, that the transition $k_{B} \to 0$ from the equilibrium statistical physics to the 
phenomenological thermodynamics is understood by identifying 
$S (E, V)$ with the entropy of the limiting thermodynamic system 
at the equilibrium state corresponding to the internal energy $E$ and the macroscopic volume $V$. 
Since on the (Lagrangian) submanifold $\Lambda \subset 
\mathbb{R}^4 (S, V; T, p)$ of the equilibrium states we have 
$d S|_{\Lambda} = ((1/ T) d E + (p / T) d V)|_{\Lambda}$, where $p$ is the macroscopic 
pressure, and $T$ is the absolute temperature, an elementary computation yields: 
\begin{equation*}
- \frac{\partial^2 S (E, V)}{\partial E^2} = 
- \Big( \frac{\partial}{\partial E} \frac{1}{T} \Big)_{V} = 
\frac{1}{T^2} \Big( \frac{\partial T}{\partial E} \Big)_{V} = 
\frac{1}{T^2 \Big( \frac{\partial E}{\partial T} \Big)_{V}} = 
\frac{1}{T^2 c_{V} (T, V)}, 
\end{equation*}
where $c_{V} (T, V)$ is the isohoric heat capacity. The assumption that 
we actually need to be able to apply the saddlepoint method $- \partial^2 S (E, V)/ \partial E^2 > 0$ 
reduces to $c_{V} (T, V) > 0$, which is one of the two conditions 
(along with $(\partial p/ \partial V)_{T} < 0$) of stability of the 
equilibrium state $(T, V)$ with respect to quasithermodynamic fluctuations.

Now let us turn to the possibility of defining a $q$-generalization of the Legendre 
transformation having the context of the four examples described above. 
Let $q = \| q_{i, j} \|_{i, j = 1}^{2 s}$ be a $2 s \times 2 s$ matrix of formal variables 
$q_{i, j}$ satisfying 
\begin{equation} 
\label{eq:qq_assumptions}
q_{i, i} = 1, \quad 
q_{i, j} q_{j, i} = 1, 
\end{equation}  
where $i, j = 1, 2, \dots, 2 s$.
Consider (over the basefield $\mathbb{K}$) an algebra $\mathcal{A}_{q}$ defined by $2 s$ generators 
$\xi_1, \xi_2, \dots, \xi_{2 s}$ and relations 
\begin{equation} 
\label{eq:xixi_rels}
\xi_{i} \xi_{j} = q_{j, i} \xi_{j} \xi_{i}, 
\end{equation}
where $i, j = 1, 2, \dots, 2 s$ (the quantum affine space).  
This algebra is $\mathbb{Z}$-graded, $\mathcal{A}_q = 
\bigoplus_{n \in \mathbb{Z}} \mathcal{A}_{q}^{n}$, 
where $\mathcal{A}_{q}^{n}$ is formed by degree $n$ homogeneous polynomials, if $n \geqslant 0$,  
and we put $\mathcal{A}_{q}^{n}$ in case $n < 0$. 
Let $\wh{\mathcal{A}}_{q}$ be the completion of $\mathcal{A}_{q}$ with respect to the 
canonical increasing filtration $F^{\bullet} \mathcal{A}_{q}$ 
associated to this grading, $F^N \mathcal{A}_{q} := 
\bigoplus_{n \leqslant N} \mathcal{A}_{q}^{n}$.  
If we think about the first $s$ generators $\xi_1, \xi_2, \dots, \xi_s$ as momenta 
(and redenote them as $\xi_1 = p_1, \xi_2 = p_2, \dots, \xi_{s} = p_{s}$), 
and the other $s$ generators $\xi_{s + 1}, \xi_{s + 2}, \dots, \xi_{2 s}$ as coordinates 
(and redenote them as $\xi_{s + 1} = x_{1}, \xi_{s + 2} = x_{2}, \dots, \xi_{2 s} = x_{s}$), 
can we define an analogue of the Legendre transformation? 
The first two examples above imply, that one could think of this analogue as a map 
$\mathcal{L}_{q} : \wh{\mathcal{A}}_{q} \to \wh{\mathcal{A}}_{q}$. 
On the other hand, the third and the fourth examples are slightly of different nature: 
there is an extra parameter introduced in the story (the Planck constant $\hbar$ in example three, 
and the Boltzmann constant $k_{B}$ in example four). 
This suggests, that the required analogue could be a map 
$\mathcal{L}_{q}^{\mathcal{B}} : \wh{\mathcal{A}}_{q} \otimes_{\mathbb{K}} \mathcal{B} \to 
\wh{\mathcal{A}}_{q} \otimes_{\mathbb{K}} \mathcal{B}$, where $\mathcal{B}$ is another 
algebra extending the scalars, or even more general, a map 
$\wt{\mathcal{L}}_{q} : \wt{\mathcal{}A}_{q} \to \wt{\mathcal{A}}_{q}$, where 
$\wt{\mathcal{A}}_{q}$ is a filtered algebra such that the degree zero of its associated 
graded is isomorphic to $\mathcal{A}_{q}$. 
The answer suggested in the present paper corresponds (more or less) to this third possibility. 

In one of our recent papers 
\cite{RuugeVanOystaeyen}, among other things, we were interested in motivating a 
$q$-analogue of the Weyl quantization map $W$ in quantum mechanics.  
Recall, that the map $W$ can be described as follows. 
One considers two algebras, a commutative algebra $\mathcal{A}$ of polynomials 
in $x_1, x_2, \dots, x_{s}$ and $p_1, p_2, \dots, p_s$, and the Heisenberg algebra 
$\wh{\mathcal{A}}$ generated by $\wh{x}_1, \wh{x}_2, \dots, \wh{x}_s, \wh{p}_1, 
\wh{p}_2, \dots, \wh{p}_s$, and 
relations $[\wh{p}_j, \wh{x}_{k}] = h$, $[\wh{p}_j, \wh{p}_{k}] = 0$, 
$[\wh{x}_j, \wh{x}_{k}] = 0$, and 
$[x_j, h] = 0$, $[p_{k}, h] = 0$, where $[-, -]$ denotes the commutator, and 
$j, k = 1, 2, \dots, s$. 
In the ``coordinate representation'' $\wh{x}_j \mapsto x_j$ (multiplication), 
$\wh{p}_k \mapsto - \ii \hbar \partial / \partial x_k$, $j, k = 1, 2, \dots, s$, and 
the central generator $h$ corresponds to the multiplication by $- \ii \hbar$. 
The map $W$ is just a linear map $W : \mathcal{A} \to \wh{\mathcal{A}}$ 
implementing the symmetrization over the ``order of action'' of the quantized 
coordinates and momenta,  
\begin{equation*} 
W (z_{i_1} z_{i_2} \dots z_{i_n}) := 
\frac{1}{n!}
\sum_{\sigma \in S_{n}} 
\wh{z}_{i_{\sigma (1)}}
\wh{z}_{i_{\sigma (2)}}
 \dots 
\wh{z}_{i_{\sigma (n)}}, 
\end{equation*}
where $S_{n}$ is the symmetric group on $n$ symbols ($n$ is a positive integer),  
$i_1, i_2, \dots, i_{n} \in \lbrace 1, 2, \dots, 2 s \rbrace$, and 
$z_{j} := p_{j}$, 
$z_{s + j} := x_{j}$, and 
$\wh{z}_j := \wh{p}_j$, 
$\wh{z}_{s + j} := \wh{x}_{j}$, for $j = 1, 2, \dots, s$. 
 
Now, being interested in a $q$-analogue of $W : \mathcal{A} \to \wh{\mathcal{A}}$, 
it was natural to introduce the so-called ``non-commutative Planck constants''. 
Note, that a similar construction motivated by the superstring theory appears in \cite{Castro}.   
The $q$-analogue mentioned is understood as a vector space map 
$W_{q} : \mathcal{A}_{q} \to \wh{\mathcal{A}}_{q}$ between the two algebras.  
The first one is the affine quantum space $\mathcal{A}_{q}$, 
denote the generators $\xi_{1}, \xi_{2}, \dots, \xi_{2 s}$, 
$\xi_{i} \xi_{j} = q_{j, i} \xi_{j} \xi_{i}$, where 
$q_{i, j}$ are formal variables satisfying the usual assumptions 
(\ref{eq:qq_assumptions}), $i, j = 1, 2, \dots, 2 s$. 
The second algebra $\wh{\mathcal{A}}_{q}$ has generators 
$\wh{\xi}_{i}$, $i = 1, 2, \dots, 2 s$, and $\wh{h}_{j, i}$, $1 \leqslant i < j \leqslant 2 s$. 
We extend the notation $\wh{h}_{j, i}$ for any $i, j = 1, 2, \dots, 2 s$, by 
$\wh{h}_{i, i} := 0$ and $h_{i, j} := - q_{j, i}^{-1} \wh{h}_{j, i}$. 
Part of the relations is the deformation of the relations 
(\ref{eq:xixi_rels}) 
for the algebra $\mathcal{A}$, 
\begin{equation*} 
\wh{\xi}_{i} \wh{\xi}_{j} - q_{j, i} \wh{\xi}_{j} \wh{\xi}_{i} = \wh{h}_{i, j}, 
\end{equation*}
for any $i, j = 1, 2, \dots, 2 s$. 
Let us order the generators as follows: 
$\wh{\xi}_{i} \prec \wh{\xi}_{j}$, if $i < j$, 
$\wh{\xi}_{i} \prec \wh{h}_{j', i'}$, for any $i$ and any $i' < j'$, and 
$\wh{h}_{j, i} \prec \wh{h}_{j', i'}$, if $i < i'$, or if $i = i'$ and $j < j'$ 
(where $i, j, i' j' \in \lbrace 1, 2, \dots, 2 s \rbrace$). 
This algebra has a Poincar\'e-Birkhoff-Witt basis (with respect to $\prec$) 
if we impose the braidings 
\begin{equation*} 
\begin{gathered}
\wh{\xi}_{i} \wh{h}_{j', i'} = q_{j', i} q_{i', i} \wh{h}_{j', i'} \wh{\xi}_{i}, \\
\wh{h}_{j, i} \wh{h}_{j', i'} = q_{j', j} q_{j', i} \, q_{i', j} q_{i', i} \, 
\wh{h}_{j', i'} \wh{h}_{j, i},   
\end{gathered}
\end{equation*} 
for any $i, j, i', j' \in \lbrace 1, 2, \dots, 2 s \rbrace$. 
With these relations, we obtain an algebra $\wh{\mathcal{A}}_{q}$ for which 
there exists a reasonable analogue $W_{q} : \mathcal{A}_{q} \to \wh{\mathcal{A}}_{q}$ of the 
Weyl quantization map 
(the generators $\wh{\xi}_{1}, \wh{\xi}_{2}, \dots, \wh{\xi}_{2 s}$
are the ``$q$-quantized'' coordinates and momenta.)
Note, that the algebra $\wh{\mathcal{A}}_{q}$ is naturally filtered, with the filtration 
$F^{\bullet} \wh{\mathcal{A}}_{q}$ 
defined by $\wh{\xi}_i \in F^{0} \wh{\mathcal{A}}_{q}$, 
$\wh{h}_{j, i} \in F^{1} \wh{\mathcal{A}}_{q} \backslash 
F^{0} \wh{\mathcal{A}}_{q}$.

It turns out that this idea to introduce the additional generators 
with non-trivial braidings is quite useful in the construction of the 
$q$-Legendre transformation. 
In other words, the generalization we suggest corresponds to the examples 
three and four described above which involve the fast oscillating integrals 
($\hbar \to 0$, the stationary point method), or 
rapidly decaying integrals ($k_{B} \to 0$, the saddle point method). 
The plan of the attack is more or less as follows. 
We start with the Hamilton-Jacobi equation 
\begin{equation*}
\partial S^t (x)/ \partial t + 
H ( x, \partial S^t (x)/ \partial x ) = 0, 
\end{equation*}
which is a non-linear equation from classical mechanics 
($x \in \mathbb{R}^s$ are the coordinates of a system, 
$t \in \mathbb{R}$ is time, 
$H (x, p)$ is the Hamiltonian, 
$p \in \mathbb{R}^s$ are the canonically conjugate momenta corresponding to $x$, 
$S^{t} (x)$ is the action as a function of (the ending point) coordinates and time).  
As is well known, the solution of the Hamilton-Jacobi equation can be described 
in terms of a system of ordinary differential equations 
$\dot x = \partial H (x, p)/ \partial p$, $\dot p = - \partial H (x, p)/ \partial x$ 
with the initial conditions $x|_{t = 0} = \alpha$, 
$p|_{t = 0} = \partial S^0 (\alpha)/ \partial \alpha$ 
(the Hamiltonian system). 
We investigate how far can we go in generalizing this fact if 
the Hamilton-Jacobi equation is replaced with its analogue 
constructed in a certain way 
over the ``braided'' generators $\lbrace \wh{\xi}_{i} \rbrace_i$ and 
$\lbrace \wh{h}_{j, i} \rbrace_{i < j}$. 
If the Hamilton-Jacobi equation is perceived as a \emph{classical limit} $\hbar \to 0$ of 
the Schr\"odinger equation, then 
at this point one realizes that it is necessary to introduce the ``Planck constants'' 
$\wh{h}_{j, i}$ back in the equation in order to ``control the braidings'' 
between the symbols in the corresponding formulae. 
The next step is to look at the fact that the 
Legendre transformation $\wt{S}^t (p)$ of $S^t (x)$ satisfies again the Hamilton-Jacobi 
equation, 
\begin{equation*} 
\partial \wt{S}^t (p)/ \partial t + 
H ( - \partial \wt{S}^t (p)/ \partial p, p ) = 0. 
\end{equation*} 
Let $\Lambda^{t} := \lbrace (x, p) \,|\, p = \partial S^t (x)/ \partial x \rbrace$. 
Then on $\Lambda^t \subset \mathbb{R}_{x}^{s} \times \mathbb{R}_{p}^{s}$ we have 
$\wt{S}^t (p)|_{\Lambda^t} = (- p x + S^t (x))|_{\Lambda^t}$. 
It turns out that there is some problem to generalize this fact to the $q$-deformed case. 
At that point we will have to make precise, what we mean by the $q$-Legendre transformation, 
but whatever it is, it is important to point out, that the construction still involves 
the additional generators having a purpose to control the braidings in the formulae. 
It remains to make the third small step. 
Once we have the $q$-Legendre formulae, there is no need to interpret the auxiliary 
generators $\wh{h}_{j, i}$ as having a quantum mechanical origin. 
For example, if we specialize all the braidings as $q_{i, j} = 1$
(and by that everything becomes commutative), then 
the analogue of the Legendre transformation can be perceived as some 
construction involving the tensor algebra $T (\bigwedge^2 (V))$, 
where $V$ is the vector space of linear functions on 
$\mathbb{R}_{x}^{s} \times \mathbb{R}_{p}^{s}$. 
There is no more reason to view $\wh{h}_{i, j}$ as non-commutative \emph{Planck} constants, 
but one could say that they could be the non-commutative ($q$-commutative) \emph{Boltzmann} constants. 
Denote them $(\wh{k}_{B})_{i, j}$. 
This suggests that if one is interested in $q$-deforming 
the statistical physics of equilibrium states 
(or, more generally, in its $R$-braiding, where $R$ is a solution of the Yang-Baxter equation), 
then it is necessary to introduce 
the non-commutative Boltzmann constants 
$(\wh{k}_{B})_{*}$ in place of $k_{B}$.

\section{$q$-analogue of the Hamilton-Jacobi equation}

Let us start with a description of the $q$-analogue of the phase space 
and the $q$-analogue of the Poisson bracket. 
Fix a positive integer $s$ and consider a $2 s \times 2 s$ matrix $q$ of 
formal variables $q_{k, l}$ satisfying $q_{k, k} = 1$ and $q_{k, l} q_{l, k} = 1$, 
for $k, l \in \lbrace 1, 2, \dots, 2 s \rbrace$. 
Consider an algebra $\mathcal{A}_q$ 
defined by generators 
$\wh{p}_1, \wh{p}_2, \dots, \wh{p}_s; \wh{x}_1, \wh{x}_2, \dots, \wh{x}_s$, and relations 
\begin{equation*}
\wh{p}_i \wh{p}_j = q_{j, i} \wh{p}_{j} \wh{p}_{i}, \quad 
\wh{p}_i \wh{x}_{\alpha} = q_{s + \alpha, i} \wh{x}_{\alpha} \wh{p}_{i}, \quad 
\wh{x}_{\alpha} \wh{x}_{\beta} = q_{s + \beta, s + \alpha} \wh{x}_{\beta} \wh{x}_{\alpha}, 
\end{equation*}
where $i, j, \alpha, \beta = 1, 2, \dots, s$, (the $q$-affine phase space).  
We will also 
use the notation $\xi_1, \xi_2, \dots, \xi_{2 s}$ (without hats), and set  
$\xi_{k} := \wh{p}_{k}$, if $k \leqslant s$, and 
$\xi_{k} := \wh{x}_{k - s}$, if $k > s$ ($k \in \lbrace 1, 2, \dots, 2 s \rbrace$). 
Extend this algebra by adding more generators $h_{\alpha, i}$, $i, \alpha = 1, 2, \dots, s$, 
satisfying 
\begin{equation*} 
\wh{p}_{i} h_{\alpha, j} = q_{s + \alpha, i} q_{j, i} h_{\alpha, j} \wh{p}_{i}, \quad 
\wh{x}_{\alpha} h_{\beta, i} = q_{s + \beta, s + \alpha} q_{i, s + \alpha} h_{\beta, i} \wh{x}_{\alpha}, 
\end{equation*}
where $i, j, \alpha, \beta = 1, 2, \dots, s$. 
Denote the extended algebra $\wt{\mathcal{A}}_{q}$ and 
equip it with a bracket 
\begin{equation*} 
\langle -, - \rangle: \wt{\mathcal{A}}_q \times \wt{\mathcal{A}}_q \to \wt{\mathcal{A}}_q,  
\end{equation*} 
defined as follows. 
Let $\langle -, - \rangle$ be a bilinear map, such that 
\begin{equation*} 
\begin{gathered}
\langle \wh{p}_{i}, \wh{x}_{\alpha} \rangle := h_{\alpha, i}, \quad 
\langle \wh{x}_{\alpha}, \wh{p}_{i} \rangle := - q_{i, s + \alpha} \langle \wh{p}_{i}, \wh{x}_{\alpha} \rangle, \\
\langle \wh{p}_{i}, \wh{p}_{j} \rangle := 0, \quad 
\langle \wh{x}_{\alpha}, \wh{x}_{\beta} \rangle := 0, 
\end{gathered}
\end{equation*}
where $i, j, \alpha, \beta = 1, 2, \dots, s$. 
Extend it to the higher order monomials in $\wh{p}_i$ and $\wh{x}_{\alpha}$ as a $q$-biderivation, 
\begin{equation*} 
\langle 
\xi_{i_m} \dots \xi_{i_1}, 
\xi_{j_n} \dots \xi_{j_1}
\rangle := 
\sum_{\mu = 1}^{m} \sum_{\nu = 1}^{n} 
Q_{i_1, \dots, i_m}^{j_1, \dots, j_n} (\mu, \nu) 
\langle \xi_{\mu}, \xi_{\nu} \rangle 
\xi_{i_m} \dots \check{\xi}_{i_{\mu}} \dots \xi_{i_1}  
\xi_{j_n} \dots \check{\xi}_{j_{\nu}} \dots \xi_{j_1}, 
\end{equation*} 
where $i_1, \dots, i_m, j_1, \dots, j_n \in \lbrace 1, 2, \dots, 2 s \rbrace$, 
the check symbol atop $\xi_{i_{\mu}}$ and $\xi_{j_{\nu}}$ means that the 
corresponding factors are omitted, 
and the braiding factor 
$Q_{i_1, \dots, i_m}^{j_1, \dots, j_n} (\mu, \nu)$ is defined from 
\begin{equation*} 
\xi_{i_m} \dots \xi_{i_1} \,  
\xi_{j_n} \dots \xi_{j_1} = 
Q_{i_1, \dots, i_m}^{j_1, \dots, j_n} (\mu, \nu) \,  
\xi_{\mu} \xi_{\nu} \, 
\xi_{i_m} \dots \check{\xi}_{i_{\mu}} \dots \xi_{i_1} \, 
\xi_{j_n} \dots \check{\xi}_{j_{\nu}} \dots \xi_{j_1}.
\end{equation*}
Finally, extend the bracket to the monomials containing 
$h_{\alpha, i}$, $\alpha, i = 1, 2, \dots, s$, as 
\begin{equation*} 
\langle h_{\alpha, i} f, g \rangle := 
h_{\alpha, i}  \langle f, g \rangle, \quad 
\langle f, g h_{\alpha, i} \rangle := 
\langle f, g \rangle h_{\alpha, i}, 
\end{equation*}
for any monomials $f, g \in \wt{\mathcal{A}}_{q}$. 
This yields a bracket $\langle -, - \rangle$ on the algebra $\wt{\mathcal{A}}_{q}$ 
which can be regarded as a $q$-analogue of the Poisson bracket $\lbrace -, - \rbrace$ 
``not divided by the Planck constant $\hbar$''. 
Define $q_{g, f}$, for any monomials $f, g \in \wt{\mathcal{A}}_q$, by $f g = q_{g, f} g f$. 
It is straightforward to check that $\langle -, - \rangle$ satisfies the $q$-Jacobi identity, 
\begin{equation*} 
\langle F, \langle G, H \rangle \rangle +
q_{G, F} q_{H, F}
\langle G, \langle H, F \rangle \rangle +
q_{H, F} q_{H, G}
\langle H, \langle F, G \rangle \rangle = 0, 
\end{equation*}
for any monomials $F, G, H \in \wt{\mathcal{A}}_{q}$ in the generators 
$\lbrace \xi_{i} \rbrace_{i}$ and $\lbrace h_{\alpha, j} \rbrace_{\alpha, j}$.

Let us now recall, how the 
Hamiltonian system of equations emerges from the 
Hamilton-Jacobi equation in the classical mechanics with one degree of freedom ($s = 1$). 
We have an unknown function $S (x, t)$ (the action) in coordinates $x \in \mathbb{R}$ 
and time $t \in \mathbb{R}$. The Hamilton-Jacobi equation looks as follows: 
\begin{equation*} 
\partial S (x, t)/ \partial t + H (x, \partial S (x, t)/ \partial x) = 0, 
\end{equation*}  
where 
$H$ is a smooth function in coordinate $x$ and momentum $p$, 
$H : \mathbb{R}_{x} \times \mathbb{R}_{p} \to \mathbb{R}$, (the Hamiltonian function), 
the momentum $p$ corresponds to $\partial S (x, t)/ \partial x$. 
Suppose we have a smooth $\mathbb{R}$-valued function $X (y, t)$ in a $\mathbb{R}$-valued variable 
$y$ and time $t$, and evaluate $\partial S (x, t)/ \partial x$ at $x = X (y, t)$. 
Then for $P (y, t) := (\partial S (x, t)/ \partial x)|_{x = X (y, t)}$ we have: 
\begin{multline*} 
\frac{\partial P (y, t)}{\partial t} = 
\Big( \frac{\partial^2 S (x, t)}{\partial x \partial t} \Big) \Big|_{x = X (y, t)} + 
\frac{\partial^2 S (x, t)}{\partial x^2} \Big|_{x = X (y, t)} 
\frac{\partial X (y, t)}{\partial t} 
= \\ = 
\Big(- \frac{\partial}{\partial x} H \Big( x, 
\frac{\partial S (x, t)}{\partial x} \Big) \Big) \Big|_{x = X (y, t)} + 
\frac{\partial^2 S (x, t)}{\partial x^2} \Big|_{x = X (y, t)} \frac{\partial X (y, t)}{\partial t}
= \\ = 
- \Big( \frac{\partial H (x, p)}{\partial x} \Big|_{p = \partial S (x, t)/ \partial x} \Big)_{x = X (y, t)} + 
\frac{\partial^2 S (x, t)}{\partial x^2} \Big|_{x = X (y, t)} 
\times \\ \times 
\Big[
- \Big( \Big( \frac{\partial H (x, p)}{\partial p} \Big) \Big|_{p = \partial 
S (x, t)/ \partial x} \Big) \Big|_{x = X (y, t)} + 
\frac{\partial X (y, t)}{\partial t}
\Big]. 
\end{multline*} 
To satisfy this, it suffices to put 
\begin{equation*} 
\frac{\partial P (y, t)}{\partial t} = 
- \frac{\partial H (x, p)}{\partial x}\bigg|_{\substack{x = X (y, t), p = P (y, t)}}, \quad 
\frac{\partial X (y, t)}{\partial t} = 
\frac{\partial H (x, p)}{\partial p} \bigg|_{\substack{x = X (y, t), p = P (y, t)}}, 
\end{equation*}
which is just the Hamiltonian system of equations.

Turn now back to the algebra $\wt{\mathcal{A}}_{q}$. 
We have the ``braided'' coordinates $\wh{x}_1, \wh{x}_2, \dots, \wh{x}_s$, 
the ```braided'' momenta $\wh{p}_1, \wh{p}_2, \dots, \wh{p}_s$, 
and the auxiliary generators $\lbrace h_{\alpha, i} \rbrace_{\alpha, i = 1}^{s}$ 
(the ``$q$-distortion'' of the Planck constant $\hbar$). 
To generalize the Hamilton-Jacobi equation to the algebra $\wt{\mathcal{A}}_{q}$, 
one should generalize the \emph{symbolic computation} linking the classical 
Hamilton-Jacobi equation with the Hamiltonian system of equations defining the 
phase space trajectories. 
Observe, that in the classical case with $s = 1$, the derivative $\partial S (x, t)/ \partial x$ 
can be expressed in terms of the canonical  Poisson bracket $\lbrace -, - \rbrace$ as  
$\partial S (x, t)/ \partial x = \lbrace p, S (x, t) \rbrace$. 
Introduce the Planck constant $\hbar$ back into the classical equations: 
\begin{equation*} 
\begin{gathered}
\frac{\partial S (x, t)}{\partial x} = \hbar \Big\lbrace p, \frac{1}{\hbar} S (x, t) \Big\rbrace, \\
\frac{1}{\hbar} \frac{\partial S (x, t)}{\partial t} + 
\frac{1}{\hbar} H \Big( x, \hbar \Big\lbrace p, \frac{1}{\hbar} S (x, t) \Big\rbrace\Big) = 0. 
\end{gathered}
\end{equation*}
One can see, that everything can be expressed in terms of the 
rescaled functions $\wh{S} (x, t) = S (x, t)/ \hbar$, 
$\wh{H} (x, p) = H (x, p)/ \hbar$, and the Poisson bracket not 
divided by $\hbar$, $\langle f, g \rangle = \hbar \lbrace f, g \rbrace$ 
(where $f$ and $g$ are any smooth functions in $(x, p)$).   
It is natural to think that the analogue of the Hamilton Jacobi equation for 
the algebra $\wt{\mathcal{A}}_{q}$ must be of the shape 
\begin{equation} 
\label{eq:q_Hamilton_Jacobi}
\frac{\partial}{\partial t} \wh{S}^t (\wh{x}) + 
\wh{H} \big( \wh{x}, \langle \wh{p}, \wh{S}^t (\wh{x}) \rangle \big) = 0, 
\end{equation} 
but what is $\wh{S}^t (\wh{x})$ and $\wh{H} (\wh{x}, \wh{p})$ in this case? 

The basic idea can be formulated as follows. 
In the classical case (for simplicity, $s = 1$), 
we have the functions $S (x, t)$, $H (x, p)$, and $X (y, t)$, which can be expanded 
into formal power series: 
\begin{equation*} 
S (x, t) = \sum_{m, N = 0}^{\infty} \frac{t^m x^N}{m! N!} b_{N}^{(m)}, \quad 
H (x, p) = \sum_{K, L = 0}^{\infty} \frac{x^K p^L}{K! L!} T_{K, L}, \quad 
X (y, t) = \sum_{m, N = 0}^{\infty} \frac{t^m y^N}{M! N!} X_{N}^{(m)},  
\end{equation*}
where $b_{N}^{(m)}$, $T_{K, L}$, and $X_{N}^{(m)}$ are some coefficients. 
Let us do \emph{everything} in terms of power series. 
In particular, the left-hand side of the Hamilton-Jacobi equation can be formally 
expanded into such series, yielding an infinite collection of links between the mentioned coefficients. 
In the multidimensional case (i.e. any $s$), the formulae are totally similar, 
but one needs to use the standard notation with multi-indices $N = (N_1, N_2, \dots, N_s)$, 
$N! = N_1! N_2! \dots N_s!$, $x^N = x_{1}^{N_1} x_{2}^{N_2} \dots x_{s}^{N_s}$, etc. 
Informally, the approach reduces to an imperative: expand everything into power series and then make the 
coefficients non-commutative.     

To illustrate the corresponding computation, assume for simplicity that 
$q_{i, j} = 1$ and $q_{s + \alpha, s + \beta} = 1$, for $i, j, \alpha, \beta = 1, 2, \dots, s$ 
(the other $q_{i', j'}$ can be non-trivial). 
Then, for example, $\langle p_i, x_{\alpha} x_{\beta} \rangle = 
h_{\alpha, i} x_{\beta} + h_{\beta, i} x_{\alpha}$ (for any $i, \alpha, \beta = 1, 2, \dots, s$), 
since one can freely permute $x_{\alpha}$ and $x_{\beta}$. 
One obtains: 
\begin{equation*} 
\wh{S}^t (\wh{x}) = \sum_{m, N}
\frac{t^m \wh{x}^N}{m! N!} b_{N}^{(m)}, \quad 
\langle \wh{p}_i, \wh{S}^t (\wh{x})\rangle = 
\sum_{m = 0}^{\infty} \frac{t^m}{m!} 
\sum_{\mu = 1}^{s} \sum_{N \geqslant 0} 
\frac{N_{\mu} h_{\mu, i} \wh{x}^{N - 1_{\mu}}}{N!} b_{N}^{(m)}, 
\end{equation*}
where $1_{\mu} := (0, \dots, 0, 1, 0, \dots, 0)$ with $1$ standing in the $\mu$-th position, 
the notation $N \geqslant 0$ is understood as $N_{\nu} \geqslant 0$ for every $\nu = 1, 2, \dots, s$. 
The time $t$ is kept as a commutative variable, although it is possible 
to consider some braidings involving $t$ as well. 
Now, since there is a factor $N_{\mu}$ in the numerator, one can assume that for $N_{\mu}$ 
the summation starts with 1, but not 0. Changing summation index to $N' = N + 1_{\mu}$, and 
then leaving out the prime in the final expression, we obtain: 
\begin{equation*} 
\langle \wh{p}_i, \wh{S}^t (\wh{x})\rangle = 
\sum_{\mu = 1}^{s} \langle \wh{p}_{i}, \wh{x}_{\mu} \rangle 
\sum_{m = 0}^{\infty} \sum_{N \geqslant 0} \frac{t^m \wh{x}^N}{m! N!} b_{N + 1_{\mu}}^{(m)}, 
\end{equation*}
where we go back to $\langle \wh{p}_i, \wh{x}_{\mu} \rangle = h_{\mu, i}$
Similarly, 
\begin{equation*} 
\frac{\partial}{\partial t} \wh{S}^t (\wh{x}) = 
\sum_{m = 0}^{\infty} 
\sum_{N \geqslant 0} 
\frac{t^m \wh{x}^N}{m! N!} b_{N}^{(m + 1)}. 
\end{equation*}
We would like the quantity $\langle p_i, \wh{S}^t (\wh{x}) \rangle$ to ``behave like'' 
$\wh{p}_i$ with respect to the braiding coefficients. 
This implies that, for every $\mu$, $m$ and $N$,  
$\langle p_i, x_{\mu} \rangle \wh{x}^{N} b_{N + 1_{\mu}}^{(m)}$ behaves like $\wh{p}_i$. 
Since the ``braiding behaviour'' of $\langle p_i, x_{\mu} \rangle$ is like of the 
product $p_i x_{\mu} = p_i x^{1_{\mu}}$, one can conclude, that for any $m$ and $N$, 
the coefficient $b_{N}^{(m)}$ behaves like $\wh{x}^{- N}$. Making it more formal, 
impose the following rules: 
\begin{equation}
\label{eq:b_xph_braiding}
\wh{p}_{i} b_{N}^{(m)} = 
\Big( \prod_{\alpha = 1}^{s} q_{s + \alpha, i}^{- N_{\alpha}} \Big)
b_{N}^{(m)} \wh{p}_{i}, \quad 
\wh{x}_{\beta} b_{N}^{(m)} = b_{N}^{(m)} \wh{x}_{\beta}, \quad 
h_{\beta, i} b_{N}^{(m)} = 
\Big( \prod_{\alpha = 1}^{s} q_{s + \alpha, i}^{- N_{\alpha}} \Big)
b_{N}^{(m)} h_{\beta, i}, 
\end{equation}  
for any $i, \beta = 1, 2, \dots, s$, any $m \in \mathbb{Z}_{\geqslant 0}$, and 
any $N = (N_1, N_2, \dots, N_s) \geqslant 0$. 

It is possible to drop the assumption the all the $\lbrace \wh{x}_{\alpha} \rbrace_{\alpha}$ commute. 
In this case, one needs to fix the order in which to write $\wh{x}_1, \wh{x}_2, \dots, \wh{x}_s$, 
say, $\wh{x}^{N} = \wh{x}_{s}^{N_s} \dots \wh{x}_2^{N_2} \wh{x}_{1}^{N_{1}}$, 
for $N = (N_1, N_2, \dots, N_s)$. 
The modification of the braiding coefficients is straightforward, and we leave 
more general formulae for the next section. 
So far, let us keep the assumption $q_{i, j} = 1$ and $q_{s + \alpha, s + \beta} = 1$, 
$i, j, \alpha, \beta = 1, 2, \dots, s$. 

It follows, that $\wh{S}^t (\wh{x})$ behaves simply as $\wh{x}^{\bar 0}$, 
$\bar 0 := (0, 0, \dots, 0)$ (length $s$), 
i.e. is a ``braided scalar''. 
From the formula for the $q$-Hamilton-Jacobi equation 
(\ref{eq:q_Hamilton_Jacobi}) 
suggested above, one 
concludes that $\wh{H} (\wh{x}, \wh{p})$ is also a ``braided scalar'', i.e. 
the coefficients $T_{K, L}$ must behave like $\wh{x}^{- K} \wh{p}^{L}$, so 
the commutation relations must be as follows: 
\begin{equation} 
\label{eq:pxT_braiding}
\wh{p}_i T_{K, L} = \Big( \prod_{\mu = 1}^{s} q_{s + \mu, i}^{- K_{\mu}} \Big) 
T_{K, L} \wh{p}_{i}, \quad \wh{x}_{\alpha} T_{K, L} =
\Big( \prod_{j = 1}^{s} q_{j, s + \alpha}^{L_{j}} \Big) T_{K, L}
\wh{x}_{\alpha},
\end{equation}
for any multi-indices $K$ and $L$, and any $i$ and $\alpha$. 
There are also commutation relations with $b_{N}^{(m)}$ and between 
different $T_{K, L}$, $T_{K', L'}$ themselves, and between different $b_{N}^{(m)}$, $b_{N'}^{(m')}$. 
Suppose we have four integer multi-indices $K$, $L$, $M$, and $N$ (of length $s$ each). 
One can look at $(\wh{x}^{K} \wh{p}^{L}) (\wh{x}^{M} \wh{p}^{N})$ and at 
$(\wh{x}^{M} \wh{p}^{N}) (\wh{x}^{K} \wh{p}^{L})$. Set 
\begin{equation} 
\label{eq:Q_KLMN}
Q_{(M, N), (K, L)} := 
\Big( \prod_{i, \alpha = 1}^{s} 
q_{i, s + \alpha}^{- N_{i} K_{\alpha}} \Big) 
\Big[ 
\prod_{j, \beta = 1}^{s} 
q_{s + \beta, j}^{L_{j} M_{\beta}}
\Big]. 
\end{equation}  
The braiding relations can be expressed as follows: 
\begin{equation*} 
\wh{p}_i T_{K, L} = Q_{(-K, -L), (\bar 0, 1_i)} T_{K, L} \wh{p}_i, \quad 
\wh{x}_{\alpha} T_{K, L} = Q_{(-K, -L), (1_{\alpha}, \bar 0)} T_{K, L} \wh{x}_{\alpha}, 
\end{equation*}
(for any $K$, $L$, $i$, $\alpha$), 
where $\bar 0 = (0, 0, \dots, 0)$ is the multi-index (of length $s$) with all zero entries. 
The other relations are as follows: 
\begin{equation} 
\label{eq:hTb_braiding}
\begin{gathered}
h_{\alpha, i} T_{K, L} = Q_{(-K, -L), (1_{\alpha, 1_i})} T_{K, L} h_{\alpha_i}, \quad 
T_{K, L} T_{K', L'} = Q_{(-K', -L'), (-K, -L)} T_{K', L'} T_{K, L}, \\ 
b_{N}^{(m)} T_{K, L} = Q_{(K, L), (-N, \bar 0)} T_{K, L} b_{N}^{(m)}, \quad 
b_{N}^{(m)} b_{N'}^{(m')} = b_{N'}^{(m')} b_{N}^{(m)}, \\
h_{\alpha, i} b_{N}^{(m)} = Q_{(-N, \bar 0), (1_{\alpha}, 1_{i})} b_{N}^{(m)} h_{\alpha, i}, \quad 
h_{\alpha, i} h_{\beta, j} = Q_{(1_\beta, 1_j), (1_{\alpha}, 1_{i})} h_{\beta, j} h_{\alpha, i}, 
\end{gathered}
\end{equation} 
for any multi-indices $K, L, K', L', N, N'$, any 
$i, j, \alpha, \beta \in \lbrace 1, 2, \dots, s \rbrace$, and any 
$m \in \mathbb{Z}_{\geqslant 0}$. 
Substitute now all the expansions into the $q$-Hamilton-Jacobi equation 
(\ref{eq:q_Hamilton_Jacobi}): 
\begin{equation} 
\label{eq:q_Hamilton_Jacobi_series}
\sum_{m, N} \frac{t^m \wh{x}^N}{m! N!} b_{N}^{(m + 1)} + 
\sum_{K, L} \frac{\wh{x}^K}{K! L!} \Big( \prod_{i = 1}^{s} \Big[
\sum_{\mu = 1}^{s} \langle \wh{p}_{i}, \wh{x}_{\mu} \rangle 
\sum_{m, N} \frac{t^m \wh{x}^N}{m! N!} b_{N + 1_{\mu}}^{(m)} \Big]^{L_{i}} \Big) T_{K, L} = 0. 
\end{equation}
Equating the terms corresponding to every $t^m \wh{x}^{N}$, it is straightforward 
to express the collection $\lbrace b_{N}^{(m + 1)} \rbrace_{N}$ in terms of 
$\lbrace b_{N}^{(m')} \rbrace_{N}$, $m' \leqslant m$, and the coefficients $h_{\alpha, i}$ and $T_{K, L}$. 
One obtains: 
\begin{multline*} 
b_{N}^{(m + 1)} = - m! N! \wh{x}^{-N} 
\sum_{K, L \in \mathbb{Z}_{\geqslant 0}^{s}} \frac{\wh{x}^K}{K! L!} 
\Big\lbrace 
\prod_{i = 1}^{s} 
\Big[
\sum_{\mu_{1}^{(i)}, \dots, \mu_{L_i}^{(i)} = 1}^{s} 
\sum_{\substack{m_{1}^{(i)}, \dots, m_{L_{i}}^{(i)} = 0,\\
\sum_{i' = 1}^{s} \sum_{r' = 1}^{L_{i'}} m_{r'}^{(i')} = m
}}^{\infty} 
\times \\ \times 
\sum_{\substack{N_{1}^{(i)}, \dots, N_{L_i}^{(i)} \in \mathbb{Z}_{\geqslant 0}^{s},\\
K + \sum_{i' = 1}^{s} \sum_{r' = 1}^{L_{i}} N_{r'}^{(i)} = N
}} 
\prod_{r = 1}^{L_{i}} 
\Big( 
h_{\mu_{r}^{(i)}, i}
\frac{\wh{x}^{N_{r}^{(i)}}}{m_{r}^{(i)}! N_{r}^{(i)}!}
\Big)
\Big]
\Big\rbrace T_{K, L}, 
\end{multline*} 
for every $N \in \mathbb{Z}_{\geqslant 0}^{s}$ and $m \in \mathbb{Z}_{\geqslant 0}$, 
and the components of the multi-index $L$ are denoted $L_i$ ($i = 1, 2, \dots, s$),   
$L = (L_1, L_2, \dots, L_s)$. 
If one brings all the factors $\wh{x}^{K}$ and 
$\wh{x}_{r}^{N_{r}^{(i)}}$, $r = 1, 2, \dots, L_i$, $i = 1, 2, \dots, s$,  
in the right-hand side to the left most position 
using the braiding relations, then one obtains 
$\wh{x}^{K + \sum_{i = 1}^{s} \sum_{r = 1}^{L_i} N_{r}^{(i)}}$, which is cancelled out 
by the factor $\wh{x}^{-N}$. 
Therefore, the formula for $b_{N}^{(m + 1)}$, in fact, 
does not contain the products of the variables 
$\wh{x}_{\alpha}$, $\alpha = 1, 2, \dots, s$.  

\begin{prop} 
Let the braiding relations 
for $b_{N}^{(m)}$ described by 
(\ref{eq:b_xph_braiding}), 
(\ref{eq:pxT_braiding}), 
(\ref{eq:hTb_braiding})
hold for $m = 0$. 
Then they are valid for any $m \in \mathbb{Z}_{\geqslant 0}$, 
if $\lbrace b_{N}^{(m)} \rbrace_{N, m}$ satisfy (\ref{eq:q_Hamilton_Jacobi_series}). 
\end{prop}

\noindent\emph{Proof.}
This follows straightforward from the 
formal analogue (\ref{eq:q_Hamilton_Jacobi_series}) of the Hamilton-Jacobi equation 
written in the form similar to $S_t = - H (x, S_x)$.  
\qed

Therefore, we arrive at the following construction. 
Let $s$ be a fixed positive integer 
(the analogue of the number of degrees of freedom 
in classical mechanics), and 
$q = \|q_{k, l}\|_{k, l = 1}^{2 s}$ be a matrix 
of formal variables. 
One considers an algebra $\mathcal{X}_{q}$ 
with an \emph{infinite} number of generators 
\begin{equation*} 
h_{\alpha, i}, \quad 
b_{N}^{(0)}, \quad 
T_{K, L}, 
\end{equation*}
where $\alpha$ and $i$ vary over 
$[s] := \lbrace 1, 2, \dots, s \rbrace$, 
$N$, $K$, and $L$ vary over $\mathbb{Z}_{\geqslant 0}^{s}$ 
(multi-indices of length $s$). 
The relations are just the braidings described above 
(see (\ref{eq:b_xph_braiding}), (\ref{eq:hTb_braiding})), 
i.e. the algebra $\mathcal{X}_{q}$ is a special case of 
an infinite-dimensional affine quantum space.  
One is interested in considering the collections 
$( f_{N} )_{N \in \mathbb{Z}_{\geqslant 0}^{s}}$ 
of elements $f_{N} \in \mathcal{X}_{q}$ satisfying 
the same braiding relations as $(b_{N}^{(0)})_{N \in 
\mathbb{Z}_{\geqslant 0}^{s}}$. 
Denote the set of all such collections as 
$\mathcal{S}_{q} (\wh{x})$. 
Then the $q$-analogue of a solution of the Hamilton-Jacobi 
equation can be perceived as an infinite 
sequence $S^{(0)}, S^{(1)}, S^{(2)}, \dots$ of 
elements $S^{(m)} \in \mathcal{S}_{q} (\wh{x})$, 
$m \geqslant 0$.  

Now let us investigate how to generalize the 
symbolic computation which derives the Hamiltonian 
system of equations $\dot p = - H_x$, $\dot x = H_p$, 
from the Hamilton-Jacobi equation $S_t + H (x, S_x) = 0$. 
We still have another formal power series  
for $X (y, t)$ that was not used. 
Assuming $s$ is generic, let 
\begin{equation*}
X (y, t) = 
\sum_{m, N} \frac{t^m y^N}{m! N!} X_{N}^{(m)},
\end{equation*} 
where $y = (y_1, y_2, \dots, y_s)$ is a collection of 
new non-commutative variables, 
$X (y, t)= (X_1 (y, t), X_2 (y, t), \dots, X_s (y, t))$, 
$X_{N}^{(m)} (y, t)= (X_{1, N}^{(m)} (y, t), X_{2, N}^{(m)} (y, t), \dots, X_{s, N}^{(m)} (y, t))$, 
$m \in \mathbb{Z}_{\geqslant 0}$, 
$N \in \mathbb{Z}_{\geqslant 0}^{s}$, 
and $t \in \mathbb{R}$. 
If one requires that $X_{\alpha} (y, t)$ behaves like 
$\wh{x}_{\alpha}$, then it means that 
$y^{N} X_{\alpha, N}^{(m)}$ behaves like $\wh{x}_{\alpha}$, 
for every $m, N$. For example, 
$p_i (y^N X_{\alpha, N}^{(m)}) = 
q_{s + \alpha, i} (y^N X_{\alpha, N}^{(m)}) p_i$, 
$i \in [s]$, etc. 
Suppose we have $\wh{S}^t (\wh{x})$ corresponding to a 
solution $((b_{N}^{(m)})_N)_{m = 0}^{\infty}$ of the 
formal $q$-Hamilton-Jacobi equation 
(\ref{eq:q_Hamilton_Jacobi_series}). 
One would like to perceive  
$P_{i} (y, t) := \langle \wh{p}_{i}, \wh{S}^t (\wh{x}) 
\rangle|_{\wh{x} = X (y, t)}$, $i = 1, 2, \dots, s$, as an analogue of a value of the 
momentum $p_{i}$ canonical conjugate to $x_{i}$, 
while the analogue of a value of $x_{i}$ is $X_{i} (y, t)$.  
We have already obtained the expressions for $\wh{S}^t (\wh{x})$ 
in a shape of a power series in $t$ and $\wh{x}$. Observe that 
the braiding relations 
(\ref{eq:b_xph_braiding}) 
between $b_{N}^{(m)}$ and $\wh{x}_{\alpha}$ imply that 
one can write $\wh{S}^{t} (\wh{x})$ both as 
$\wh{S}^{t} (\wh{x}) = \sum_{m, N} (t^m \wh{x}^{N}/ m! N!) b_{N}^{(m)}$, and as 
$\wh{S}^{t} (\wh{x}) = \sum_{m, N} b_{N}^{(m)} (t^m \wh{x}^{N}/ m! N!)$, 
and that $\langle \wh{p}_i, \wh{S}^t (\wh{x}) \rangle = - \langle \wh{S}^t (\wh{x}), \wh{p}_i \rangle$, 
$i \in [s]$. 
Similar remarks can be made for $\wh{H} (\wh{x}, \wh{p}) = 
\sum_{K, L} (\wh{x}^K \wh{p}^L/ K! L!) T_{K, L}$. Applying the bracket $\langle \wh{p}_i, - \rangle$, 
one obtains: 
\begin{equation*} 
\langle \wh{p}_i, \wh{H} (\wh{x}, \wh{p}) \rangle = 
\sum_{\mu = 1}^{s} \langle \wh{p}_i, \wh{x}_{\mu} \rangle 
\sum_{K, L \in \mathbb{Z}_{\geqslant 0}^{s}} 
\frac{\wh{x}^{K} \wh{p}^{L}}{K! L!} T_{K + 1_{\mu}, L}, 
\end{equation*}
for each $i \in [s] = \lbrace 1, 2, \dots, s \rbrace$, 
and we keep $\langle \wh{p}_i, \wh{x}_{\mu} \rangle$ instead of writing $h_{\mu, i}$. 
One can compute similarly $\langle \wh{x}_{\alpha}, \wh{H} (\wh{x}, \wh{p}) \rangle$, 
but now it is more convenient to use the expression 
$\wh{H} (\wh{x}, \wh{p}) = 
\sum_{K, L}  T_{K, L} \wh{x}^K \wh{p}^L/ (K! L!)$, and to apply the bracket with $\wh{x}_{\alpha}$ 
from the right, $\langle \wh{x}_{\alpha}, \wh{H} (\wh{x}, \wh{p}) \rangle = 
- \langle \wh{H} (\wh{x}, \wh{p}), \wh{x}_{\alpha}\rangle$, 
\begin{equation*} 
\langle \wh{x}_{\alpha}, \wh{H} (\wh{x}, \wh{p}) \rangle = 
- \sum_{j = 1}^{s} \sum_{K, L \in \mathbb{Z}_{\geqslant 0}^{s}} T_{K, L + 1_{j}} 
\frac{\wh{x}^K \wh{p}^L}{K! L!} \langle \wh{p}_j, \wh{x}_{\alpha} \rangle, 
\end{equation*}
for each $\alpha \in [s]$, $\langle \wh{p}_j, \wh{x}_{\alpha} \rangle = h_{\alpha, j}$. 
In what follows, in case it does not lead to confusion, 
it is convenient to drop the summation over $j \in [s]$ from the notation, 
and write (in analogy with the matrix multiplication) simply 
$\langle x_{\alpha}, \wh{H} (\wh{x}, \wh{p}) \rangle = 
- \sum_{K, L} T_{K, L + 1} (\wh{x}^K \wh{p}^L/ K! L!) \langle \wh{p}, \wh{x}_{\alpha} \rangle$.

Now, if we look at the $q$-analogue of the Hamilton-Jacobi equation 
(\ref{eq:q_Hamilton_Jacobi_series}), 
applying $\langle \wh{p}_i, - \rangle$ to the left and the right-hand sides, yields: 
\begin{multline*} 
\langle \wh{p}_i, \wh{x} \rangle 
\sum_{m, N} \frac{t^m \wh{x}^N}{m! N!} b_{N + 1}^{(m)} + 
\langle \wh{p}_i, \wh{x} \rangle \sum_{K, L} 
\frac{\wh{x}^K}{K! L!} 
\langle \wh{p}, \wh{S}^t (\wh{x}) \rangle^{L} T_{K + 1, L} 
- \\ - 
\sum_{K, L} T_{K, L + 1} \frac{\wh{x}^K}{K! L!} 
\langle \wh{p}, \wh{S}^t (\wh{x}) \rangle^{L} 
\langle \langle \wh{p}, \wh{S}^t (\wh{x}) \rangle, \wh{p}_i \rangle = 0. 
\end{multline*}
We have a power series $X (y, t)$, which corresponds to the coordinates $\wh{x}$. 
Consider the series for $\langle \wh{p}, \wh{S}^t (\wh{x}) \rangle$ and then 
substitute in it $\wh{x} \to X (y, t)$. This defines a quantity 
$P (y, t) = (P_1 (y, t), P_2 (y, t), \dots, P_s (y, t))$, 
\begin{equation*}
P (y, t) := \langle \wh{p}, \wh{S}^t (\wh{x}) \rangle|_{x \to X (y, t)}
\end{equation*}
The aim is to link 
the derivatives $\partial X (y, t)/ \partial t$ and $\partial P (y, t)/ \partial t$ with 
the brackets $\langle \wh{x}, \wh{H} (\wh{x}, \wh{p}) \rangle$ and 
$\langle \wh{p}, \wh{H} (\wh{x}, \wh{p}) \rangle$. 
A straightforward computation yields: 
\begin{equation*} 
\frac{\partial P_i (y, t)}{\partial t} = 
\Big\langle \wh{p}_i, \frac{\partial \wh{S}^t (\wh{x})}{\partial t} \Big\rangle 
\Big|_{\wh{x} \to X (y, t)} + 
\sum_{\alpha, \beta \in [s]} 
\frac{\partial X_{\beta} (y, t)}{\partial t} 
\sum_{m, N} \frac{t^m X (y, t)^N}{m! N!} 
b_{N + 1_{\beta} + 1_{\alpha}}^{(m)} \langle \wh{x}_{\alpha}, \wh{p}_i \rangle. 
\end{equation*}
The first term in the right-hand side can be transformed using the result of the 
computation of $\langle \wh{p}_i, - \rangle$ applied to the $q$-Hamilton-Jacobi equation. 
It can be expressed as follows: 
\begin{multline*} 
\Big\langle \wh{p}_i, \frac{\partial \wh{S}^t (\wh{x})}{\partial t} \Big\rangle 
\Big|_{\wh{x} \to X (y, t)} = 
- (\langle \wh{p}, \wh{H} (\wh{x}, \wh{p}) 
\rangle|_{\wh{p} \to \langle \wh{p}, \wh{S}^t (\wh{x}) \rangle} )|_{\wh{x} \to X (y, t)} 
+ \\ + 
\sum_{j, \alpha, \beta \in [s]} 
\sum_{K, L} T_{K, L + 1_{j}} \frac{\wh{x}^K}{K! L!} \langle \wh{p}, \wh{S}^t (\wh{x}) \rangle^L 
\langle \wh{p}_j, \wh{x}_{\beta} \rangle 
\sum_{m, N} 
\frac{t^m \wh{x}^N}{m! N!} b_{N + 1_{\beta} + 1_{\alpha}}^{(m)} 
\langle \wh{x}_{\alpha}, \wh{p}_i \rangle. 
\end{multline*}
Substituting this expression into the previous expression 
for $\partial P_i (y, t)/ \partial t$, we obtain: 
\begin{multline*} 
\frac{\partial P_i (y, t)}{\partial t} = 
- \langle \wh{p}_i, \wh{H} (\wh{x}, \wh{p}) \rangle|_{\wh{x} \to X (y, t), \wh{p} \to P (y, t)} 
+ \\ + 
\sum_{\beta} \Big\lbrace - \langle \wh{H} (\wh{x}, \wh{p}), \wh{x}_{\beta} 
\rangle|_{\wh{x} \to X (y, t), \wh{p} \to P (y, t)} + 
\frac{\partial X_{\beta} (y, t)}{\partial t} \Big\rbrace 
\sum_{\alpha, m, N} 
\frac{t^m \wh{x}^N}{m! N!} b_{N + 1_{\beta} + 1_{\alpha}}^{(m)} 
\langle \wh{x}_{\alpha}, \wh{p}_i \rangle. 
\end{multline*}
Therefore, one can see, that it suffices to put 
\begin{equation*}
\begin{aligned} 
\frac{\partial P_i (y, t)}{\partial t} &= 
- \langle \wh{p}_i, \wh{H} (\wh{x}, \wh{p}) \rangle
\Big|_{x \to X (y, t), p \to P (y, t)}, \\
\frac{\partial X_{\alpha} (y, t)}{\partial t} &= 
- \langle \wh{x}_{\alpha}, \wh{H} (\wh{x}, \wh{p})\rangle 
\Big|_{x \to X (y, t), p \to P (y, t)}, 
\end{aligned}
\end{equation*}
where $P_i (y, t) = \langle \wh{p}_i, \wh{S}^t (\wh{x}) \rangle|_{\wh{x} \to X (y, t)}$, 
$i, \alpha \in [s] = \lbrace 1, 2, \dots, s \rbrace$. 
Recall, that these equations correspond to the $q$-Hamilton-Jacobi equation ``in 
$\wh{x}$-representation'', 
\begin{equation*} 
\frac{\partial \wh{S}^t (\wh{x})}{\partial t} + 
\wh{H} (\wh{x}, \wh{p})\big|_{\wh{p} \to \langle \wh{p}, \wh{S}^t (\wh{x}) \rangle} = 0. 
\end{equation*} 
This shape of equations is quite natural to expect, since 
for the usual Hamiltonian system expressed in terms of the Poisson bracket 
not divided by the Planck constant $\hbar$, we have 
$\dot p = - \lbrace p, H (x, p) \rbrace = - \hbar \lbrace p, \hbar^{-1} H (x, p) \rbrace$, and 
$\dot x = - \lbrace x, H (x, p) \rbrace = - \hbar \lbrace x, \hbar^{-1} H (x, p) \rbrace$. 

So far there was nothing assumed about the nature of $y = (y_1, y_2, \dots, y_s)$ 
in the formal power series $X (y, t) = \sum_{m, N} (t^m y^N/ (m! N!)) X_{N}^{(m)}$. 
We also have 
\begin{equation*} 
P_{i} (y, t) = - \sum_{\mu, m, N} \frac{t^m}{m! N!} 
\Big( 
\sum_{k, M} \frac{t^{k} y^{M}}{k! M!} X_{M}^{(k)} 
\Big)^{N} b_{N + 1_{\mu}} \langle x_{\mu}, p_{i} \rangle, 
\end{equation*}
for every $i \in [s]$. 
One can say, that $y^N X_{N}^{(m)}$ ``behaves like'' (i.e. has the same braiding coefficients) 
as $\wh{x}$, for every $m \in \mathbb{Z}_{\geqslant 0}$ and $N \in \mathbb{Z}_{\geqslant 0}^{s}$. 
Let us assume, that $y$ behaves like $\wh{x}$. In classical mechanics, this corresponds 
to considering the Hamiltonian system of equations with initial conditions 
on a Lagrangian manifold $\Lambda$, $\dim \Lambda = s$, 
which projects diffeomorphically on the configuration space 
$\mathbb{R}_{x}^{s}$ along the momenta space $\mathbb{R}_{p}^{s}$, and 
taking the configuration space coordinates as the local coordinates. 
Then one can perceive $P (y, t)$ as a formal power series 
\begin{equation*}
P (y, t) = \sum_{n, M} \frac{t^n y^M}{n! M!} P_{M}^{(n)},
\end{equation*}
where $y^{M} P_{M}^{(n)}$ behaves like $\wh{p}$ 
for every $n \in \mathbb{Z}_{\geqslant 0}$ and $M \in \mathbb{Z}_{\geqslant 0}^{s}$.  
Consider the $q$-Hamiltonian system with these $X (y, t)$ and $P (y, t)$: 
\begin{equation} 
\label{eq:q_Hamiltonian_XPyt}
\frac{\partial P (y, t)}{\partial t} = \langle \wh{H}, \wh{p} \rangle, \quad 
\frac{\partial X (y, t)}{\partial t} = \langle \wh{H}, \wh{x} \rangle, 
\end{equation}
where the right-hand sides are evaluated at $\wh{x} \to X (y, t)$, $\wh{p} \to P (y, t)$.
Then one obtains: 
\begin{multline*} 
\sum_{n, M} \frac{t^n y^M}{n! M!} P_{M}^{(n + 1)} = 
- \sum_{\mu = 1}^{s} \langle \wh{p}, \wh{x}_{\mu} \rangle
\sum_{K, L} \Big( 
\prod_{\alpha = 1}^{s} \Big[ 
\sum_{m^{(\alpha)}, N^{(\alpha)}} 
\frac{t^{m^{(\alpha)}} y^{(N^{(\alpha)})}}{m^{(\alpha)}! N^{(\alpha)}!}
X_{\alpha, N^{(\alpha)}}^{(m^{(\alpha)})}
\Big]^{K_{\alpha}}
\Big) 
\times \\ \times 
\Big( 
\prod_{i = 1}^{s} \Big[ 
\sum_{n^{(i)}, M^{(i)}} 
\frac{t^{n^{(i)}} y^{(M^{i})}}{n^{(i)}! M^{(i)}!}
P_{i, M^{(i)}}^{(n^{(i)})}
\Big]^{L_{i}}
\Big) T_{K + 1_{\mu}, L}, 
\end{multline*}
where we put $P_{M}^{(n)} = (P_{1, M}^{(n)}, P_{2, M}^{(n)}, \dots, P_{s, M}^{(n)})$ and 
$X_{N}^{(m)} = (X_{1, M}^{(m)}, X_{2, M}^{(n)}, \dots, X_{s, M}^{(n)})$ to denote the 
corresponding components. There is a similar expression corresponding to the 
second part of the $q$-Hamiltonian system 
(\ref{eq:q_Hamiltonian_XPyt}). From the latter equation, 
one can extract the coefficient $P_{M}^{(n + 1)}$ (any $n \in \mathbb{Z}_{\geqslant 0}$ and 
any $M \in \mathbb{Z}_{\geqslant 0}^{s}$) in the left-hand side 
by restricting the summations in the right-hand side to obtain the same powers by $t$ and $y$.  
It is convenient to introduce,   
for every $m, n \in \mathbb{Z}_{\geqslant 0}$ and each $M, N \in \mathbb{Z}_{\geqslant 0}^{s}$, 
the quantities  
\begin{equation*} 
\bar P_{M}^{(n)} := y^{M} P_{M}^{(n)}, \quad 
\bar X_{N}^{(m)} := y^{N} X_{N}^{(m)}, 
\end{equation*} 
and to assemble them into generating functions 
\begin{equation*} 
\bar P (t, \varepsilon) := \sum_{n, M} \frac{t^n \varepsilon^{M}}{n! M!} \bar P_{M}^{(n)}, \quad 
\bar X (t, \varepsilon) := \sum_{m, N} \frac{t^m \varepsilon^{N}}{m! N!} \bar X_{N}^{(m)}, 
\end{equation*}
where $\varepsilon = (\varepsilon_1, \varepsilon_2, \dots, \varepsilon_s)$ is a vector 
of formal \emph{commutative} variables (i.e. every $\varepsilon_i$, $i \in [s]$,  
commutes with every other symbol in the expressions, just like the time $t$). 
Denote the components of $\bar P (t, \varepsilon)$ 
as $\bar P_i (t, \varepsilon)$, $i \in [s]$, and 
the components of $\bar X (t, \varepsilon)$ as $\bar X_{\alpha} (t, \varepsilon)$, $\alpha \in [s]$.  
We have the commutation (braiding) relations:  
\begin{equation} 
\label{eq:XPtepsilon_braiding}
\begin{gathered} 
\bar P_{i} (t, \varepsilon) \bar X_{\alpha} (t', \varepsilon') = 
Q_{(1_{\alpha}, \bar 0), (\bar 0, 1_{i})} \bar X_{\alpha} (t', \varepsilon') \bar P_{i} (t, \varepsilon), \\
\bar X_{\alpha} (t, \varepsilon) \bar X_{\beta} (t', \varepsilon') = 
Q_{(1_{\beta}, \bar 0), (1_{\alpha}, \bar 0)}
\bar X_{\beta} (t', \varepsilon') \bar X_{\alpha} (t, \varepsilon), \\
\bar P_{i} (t, \varepsilon) \bar P_{j} (t', \varepsilon') = 
Q_{(\bar 0, 1_{j}), (\bar 0, 1_{i})}
\bar P_{j} (t', \varepsilon') \bar P_{i} (t, \varepsilon), \\
\bar P_{i} (t, \varepsilon) h_{\beta, j} = 
Q_{(1_{\beta}, 1_{j}), (\bar 0, 1_{i})} h_{\beta, j} \bar P_{i} (t, \varepsilon), \quad 
\bar X_{\alpha} (t, \varepsilon) h_{\beta, j} = 
Q_{(1_{\beta}, 1_{j}), (1_{\alpha}, \bar 0)} h_{\beta, j} \bar X_{\alpha} (t, \varepsilon), \\
\bar P_{i} (t, \varepsilon) T_{K, L} = 
Q_{(K, L), (\bar 0, 1_{i})} T_{K, L} \bar P_{i} (t, \varepsilon), \quad 
\bar X_{\alpha} (t, \varepsilon) T_{K, L} = 
Q_{(K, L), (1_{\alpha}, \bar 0)} T_{K, L} \bar X_{\alpha} (t, \varepsilon), 
\end{gathered}
\end{equation}
where $i, j, \alpha, \beta \in [s]$, $K, L \in \mathbb{Z}_{\geqslant 0}^{s}$, and 
the notation $Q_{(K, L), (M, N)}$ for any multi-indices $K, L, M, N$ is introduced in 
(\ref{eq:Q_KLMN}).  
In particular if we specialize $t$ to zero, one recovers the relations involving the generators 
\begin{equation*} 
\bar P_{i, M}^{(0)}, \quad 
\bar X_{\alpha, N}^{(0)}, \quad 
h_{\alpha, i}, \quad T_{K, L}, 
\end{equation*}
where $i$ and $\alpha$ vary over $[s]$, and $K$, $L$, $M$, $N$ vary over $\mathbb{Z}_{\geqslant 0}^{s}$. 
Denote the corresponding algebra (an affine quantum space) $\mathcal{Z}_{q}$. 
One can consider the Hamiltonian system over this algebra: 
\begin{equation} 
\label{eq:q_Hamiltonian_Zq}
\frac{\partial \bar P_i (t, \varepsilon)}{\partial t} = \langle \wh{H} (\wh{x}, \wh{p}), \wh{p}_i \rangle, \quad 
\frac{\partial \bar X_{\alpha} (t, \varepsilon)}{\partial t} = 
\langle \wh{H} (\wh{x}, \wh{p}), \wh{x}_{\alpha} \rangle, 
\end{equation}
where $i, \alpha \in [s]$, and the right hand sides are evaluated at 
$\wh{x} \to \bar X (t, \varepsilon)$, $\wh{p} \to \bar P (t, \varepsilon)$. 
\begin{prop} 
If the coefficients $\bar P_{i, M}^{(n)} \in \mathcal{Z}_{q}$ and 
$\bar X_{\alpha, N}^{(m)} \in \mathcal{Z}_{q}$, $i, \alpha = 1, 2, \dots, s$, 
$n, m \in \mathbb{Z}_{\geqslant 0}$, 
$M, N \in \mathbb{Z}_{\geqslant 0}^{s}$, 
are defined by the $q$-Hamiltonian system 
(\ref{eq:q_Hamiltonian_Zq}), 
then the braiding relations 
(\ref{eq:XPtepsilon_braiding}) 
hold for any 
$t$, $t'$, $\varepsilon$, $\varepsilon'$. 
\end{prop}
\noindent
\emph{Proof.} 
Straightforward. 
\qed

Intuitively, what happens is that if one is interested in the non-commutative 
($q$-commutative) analogue of the Hamiltonian system in the classical mechanics 
$\dot p = - H_{x} (x, p)$, $\dot x = H_{p} (x, p)$, 
one needs to \emph{``blow up''} every phase space point 
$(x, p) \in \mathbb{R}_{x}^{s} \times \mathbb{R}_{p}^{s}$ into an \emph{infinite} 
collection of generators, 
\begin{equation*} 
p_{i} \to \lbrace \bar P_{i, M}^{(0)} \rbrace_{M \in \mathbb{Z}_{\geqslant 0}}, \quad 
x_{\alpha} \to \lbrace \bar X_{\alpha, N}^{(0)} \rbrace_{M \in \mathbb{Z}_{\geqslant 0}}, 
\end{equation*} 
where $i, \alpha \in [s]$. 
The algebra of classical polynomial observables is replaced by the algebra 
$\mathcal{Z}_{q}$, involving additional generators $h_{\alpha, i}$, 
$i, \alpha \in [s]$, and $T_{K, L}$, $K, L \in \mathbb{Z}_{\geqslant 0}^{s}$.

\section{Noncommutative Legendre transformation}

One can try to use a similar trick with multiplying and dividing over $\hbar$ 
and working in terms of Poisson brackets $\lbrace -, - \rbrace$ rather than the
derivatives, to define a $q$-analogue of the Legendre transformation \cite{Tulczyjew}. 
One runs into a problem here, though. If we look at a classical formula in 
one-dimensional case, $s = 1$, then we have:  
\begin{equation*} 
\wt{S} (p)|_{p = \partial S (x)/ \partial x} = - \frac{\partial S (x)}{\partial x} x + S (x). 
\end{equation*}  
Now, the idea is to replace $\partial S (x)/ \partial x$ with 
$\hbar \lbrace p, \hbar^{-1} S (x) \rbrace$. 
This yields: 
\begin{equation*} 
\hbar^{-1} \wt{S} (p)|_{p \to \hbar \lbrace p, \hbar^{-1} S (x) \rbrace} = 
- \frac{1}{\hbar} (\hbar \lbrace p, \hbar^{-1} S (x) \rbrace) x + 
\hbar^{-1} S (x). 
\end{equation*}
One can see that there is a ``singularity'' $1/ \hbar$ in front of the first term 
in the right-hand side as $\hbar \to 0$. 
If one perceives $S (x)$ and $\wt{S} (p)$ as the arguments of the rapidly 
decaying exponents in a saddle-point method, then the idea reminds an attempt to 
multiply two Dirac delta-functions concentrated in the same point. 

On the other hand, we still can follow the informal imperative 
``to expand everything in power series and then make the coefficients non-commutative''. 
Let us remind how to compute the corresponding coefficients in the 
classical one-dimensional case. 
One has: 
\begin{equation*} 
\wt{S} (p) = (- x p + S (x))|_{x = \wt{x} (p)}, 
\end{equation*}
where $\wt{x} (p)$ is defined as a solution of $p = \partial S (x)/ \partial x$ with respect 
to $x$, (we assume that this solution is unique for every value of the parameter $p$). 
Differentiating the identity 
$p \equiv (\partial S (x)/ \partial x)|_{x = \wt{x} (p)}$ over $p$, one obtains 
$1 = (\partial^2 S (x)/ \partial x^2)|_{x = \wt{x} (p)} \partial \wt{x} (p)/ \partial p$, 
and, therefore, 
\begin{equation*} 
\frac{\partial \wt{x} (p)}{\partial p} = u (x)|_{x = \wt{x} (p)}, \quad 
u (x) := \frac{1}{\partial^2 S (x)/ \partial x^2}. 
\end{equation*}
On the other hand, differentiating $\wt{S} (p)$ over $p$ yields 
\begin{equation*} 
\frac{\partial \wt{S} (p)}{\partial p} = 
- \frac{\partial \wt{x} (p)}{\partial p} p - \wt{x} (p) + 
\frac{\partial S (x)}{\partial x} \Big|_{x = \wt{x} (p)} \frac{\partial \wt{x} (p)}{\partial p} = - \wt{x} (p). 
\end{equation*}
For the the second derivative, one obtains 
$\partial^2 \wt{S} (p)/ \partial p^2 = - \partial \wt{x} (p)/ \partial p = - u (x)|_{x = \wt{x} (p)}$. 
The third derivative 
$\partial^3 \wt{S} (p)/ \partial p^3 = - u' (x)|_{x = \wt{x} (p)} \partial \wt{x} (p)/ \partial p = 
- (u (x) u' (x))|_{x = \wt{x} (p)}$. By induction, 
\begin{equation*} 
\frac{\partial^{m + 2} \wt{S} (p)}{\partial p^{m + 2}} = 
- \Big( \Big( u (x) \frac{\partial}{\partial x} \Big)^{m} u (x) \Big) \Big|_{x = \wt{x} (p)}, 
\end{equation*}
for $m \geqslant 0$. 
Let us assume that $\wt{x} (0) = 0$. 
Then we have 
\begin{equation*} 
\frac{\partial^{m + 2} \wt{S} (p)}{\partial p^{m + 2}} \Big|_{p = 0} = 
- \Big( \Big( u (x) \frac{\partial}{\partial x} \Big)^{m} u (x) \Big) \Big|_{x = 0}, \quad m \geqslant 0, 
\end{equation*} 
where the right-hand side can be expressed in terms of 
$(\partial^{2 + n} S (x)/ \partial x^{2 + n})|_{x = 0}$, $n \geqslant 0$. 
It is also possible to assume (without loss of generality), that 
$S (x)|_{x \to 0} = 0$ and 
$\wt{S} (p)|_{p \to 0} = 0$, so the corresponding Taylor expansions start at quadratic terms. 
In what follows, we shall generalize these formulae for the coefficients, 
or, more precisely, the just described way to link these coefficients.  

In this section we assume, for simplicity, that the formal variables 
$q_{k, l}$, $k, l = 1, 2, \dots, 2 s$ are all specialized to $1$. 
Note, that even in this case we still have a matrix $h = \| h_{k, l} \|_{k, l = 1}^{2 s}$, 
and not something like $\hbar (\delta_{k, l + s} - \delta_{k + s, l})$. 
The aim, in general, is to link two formal expansions starting at quadratic terms: 
\begin{equation*} 
\wh{S} (\wh{x}) = \sum_{N} \frac{\wh{x}^{N}}{N!} b_{N}, \quad 
\wh{\wt{S}} (\wh{p}) = \sum_{M} \frac{\wh{p}^{M}}{M!} a_{M}, 
\end{equation*} 
where $M$ and $N$ are multi-indices of length $s$, $b_{N}$ 
behaves like $\wh{x}^{- N}$, and $a_{M}$ ``behaves like'' $\wh{p}^{- M}$ 
(in the sense explained in the previous section). 
For convenience, in case it can not cause a confusion, we omit the hats in the notation 
$\wh{x} = (\wh{x}_1, \wh{x}_2, \dots, \wh{x}_s)$, $\wh{S} (\wh{x})$, 
$\wh{p} = (\wh{p}_1, \wh{p}_2, \dots, \wh{p}_s)$, and $\wh{\wt{S}} (\wh{p})$. 
Denote: 
\begin{equation*} 
\begin{gathered}
B_{k_1, k_2, \dots, k_r} (x) := \langle p_{k_r}, \dots 
\langle p_{k_2}, \langle p_{k_1}, S (x) \rangle \rangle \dots \rangle, \\
A_{\lambda_1, \lambda_2, \dots, \lambda_r} (p) := 
\langle x_{\lambda_r}, \dots 
\langle x_{\lambda_2}, \langle x_{\lambda_1}, \wt{S} (p) \rangle \rangle \dots \rangle,  
\end{gathered}
\end{equation*}
where $k_1, k_2, \dots, k_r \in [s]$ and $\lambda_1, \lambda_2, \dots, \lambda_r \in [s]$, 
and $r \in \mathbb{Z}_{> 0}$. 
Let the leading coefficients in the corresponding formal power series be denoted as 
\begin{equation*} 
B_{k_1, k_2, \dots, k_r}^{(0)} := B_{k_1, k_2, \dots, k_r} (x)|_{x \to \bar 0}, \quad 
A_{\lambda_1, \lambda_2, \dots, \lambda_r}^{(0)} := 
A_{\lambda_1, \lambda_2, \dots, \lambda_r} (p)|_{p \to \bar 0}. 
\end{equation*}
Let us mimic the computation of the coefficients for the classical Legendre transformation in this set-up. 

In case of the classical Legendre transformation we have a fact 
$S_{x x} \wt{S}_{p p} + 1 = 0$, where the left-hand side is evaluated at $x = \wt{x} (p)$. 
In the $q$-deformed case, $x$ corresponds to $\langle \wh{x}, \wt{\wh{S}} (\wh{p}) \rangle$, 
and $p$ corresponds to $\langle \wh{p}, \wh{S} (\wh{x}) \rangle$. 
Omitting the hats, we have: 
\begin{equation*} 
p_i = \langle p_i, S (x) \rangle|_{x \to \langle x, \wt{S} (p) \rangle} = 
\sum_{\mu = 1}^{s} \langle p_i, x_{\mu} \rangle \sum_{N} \frac{1}{N!} 
\Big[ 
\sum_{k = 1}^{s} \langle x, p_k \rangle \sum_{M} \frac{p^M}{M!} a_{M + 1_{k}}
\Big]^{N} b_{N + 1_{\mu}}, 
\end{equation*}
where $i \in [s]$. Apply $\langle -, x_{\omega} \rangle$, $\omega \in [s]$: 
\begin{multline} 
\label{eq:pxbNaM}
\langle p_{i}, x_{\omega} \rangle = 
\sum_{\mu, \nu, k, l = 1}^{s} 
\langle p_{i}, x_{\mu} \rangle 
\sum_{N} \frac{1}{N!} 
\Big[ 
\sum_{k' = 1}^{s} \langle x, p_{k'} \rangle \sum_{M'} \frac{p^{M'}}{{M'}!} a_{M' + 1_{k'}}
\Big]^{N} 
b_{N + 1_{\mu} + 1_{\nu}}
\times \\ \times 
\langle x_{\nu}, p_{k} \rangle 
\sum_{M} \frac{p^M}{M!} a_{M + 1_{k} + 1_{l}} 
\langle p_{l}, x_{\omega} \rangle. 
\end{multline}
Taking the leading term in $p$ in the right-hand side, one obtains: 
\begin{equation} 
\label{eq:pxba}
\langle p_{i}, x_{\omega} \rangle = 
\sum_{\mu, \nu, k, l = 1}^{s} 
\langle p_{i}, x_{\mu} \rangle 
b_{1_{\mu} + 1_{\nu}}
\langle x_{\nu}, p_{k} \rangle 
a_{1_{k} + 1_{l}} 
\langle p_{l}, x_{\omega} \rangle. 
\end{equation} 
Similarly, starting with 
$x_{\omega} = \langle x_{\omega}, \wt{S} (p) \rangle|_{p \to \langle p, S (x) \rangle }$, $\omega \in [s]$, 
applying the bracket $\langle -, p_{i} \rangle$, $i \in [s]$, and taking the leading coefficient, yields: 
\begin{equation} 
\label{eq:xpab}
\langle x_{\omega}, p_{i} \rangle = 
\sum_{k, l, \mu, \nu = 1}^{s}
\langle x_{\omega}, p_{k} \rangle a_{1_{k} + 1_{l}} 
\langle p_{l}, x_{\mu} \rangle b_{1_{\mu} + 1_{\nu}} 
\langle x_{\nu}, p_{i} \rangle.  
\end{equation}
The last two equalities 
(\ref{eq:pxba}), (\ref{eq:xpab})
can be perceived as analogues of 
$(S_{x x}|_{x \to 0}) (\wt{S}_{p p}|_{p \to 0}) + 1 = 0$ in the classical case. 
Now, let us introduce the collection of variables $u_{N} (\alpha, \beta)$ by 
\begin{equation} 
\label{eq:a_uN}
\sum_{k, l = 1}^{s} \langle x_{\alpha}, p_{k} \rangle 
\sum_{M}
\frac{p^M}{M!}
a_{M + 1_{k} + 1_{l}} 
\langle p_{l}, x_{\beta} \rangle = 
\sum_{N} \frac{1}{N!} \Big[ 
\sum_{i = 1}^{s}
\langle x, p_{i} \rangle \Big( \sum_{M} \frac{p^M}{M!} a_{M + 1_{i}} \Big)
\Big]^{N} u_{N} (\alpha, \beta), 
\end{equation}
where $M$ and $N$ are multi-indices of length $s$. 
The expression in the square brackets in the right-hand side is just 
$\langle x, \wt{S} (p) \rangle$. 
Note, that from the assumptions about the braiding behaviour of $x$, $p$, $b_{N}$, and $a_{M}$ 
($b_{N}$ behaves like $x^{- N}$ and $a_{M}$ behaves like $p^{- M}$), one should assume 
that $u_{N} (\alpha, \beta)$ behaves like $x^{-N + 1_{\alpha} + 1_{\beta}}$, 
for any $\alpha, \beta \in [s]$ and any $N \in \mathbb{Z}_{\geqslant 0}^{s}$.  
Furthermore, using the two previous identities, one has: 
\begin{equation*} 
\begin{aligned}
\sum_{N} \frac{x^N}{N!} u_{N} (\alpha, \beta) &= 
\sum_{k, l, \mu, \nu = 1}^{s}
\langle x_{\alpha}, p_{k} \rangle a_{1_{k} + 1_{l}} \langle p_{l}, x_{\mu} \rangle b_{1_{\mu} + 1_{\nu}} 
 \sum_{N} \frac{x^N}{N!} u_{N} (\nu , \beta), \\
\sum_{N} \frac{x^N}{N!} u_{N} (\alpha, \beta) &= 
\sum_{\mu, \nu, k, l = 1}^{s}  
\sum_{N} \frac{x^N}{N!} u_{N} (\alpha, \mu) 
b_{1_{\mu} + 1_{\nu}}
\langle x_{\nu}, p_{k} \rangle
a_{1_{k} + 1_{l}}
\langle p_{l}, x_{\beta} \rangle, 
\end{aligned}
\end{equation*}
where $\alpha, \beta \in [s]$. 
Using the first of these equalities and the formulae 
(\ref{eq:pxbNaM}) 
linking $a_{M}$ and $b_{N}$, 
one obtains: 
\begin{multline} 
\label{eq:pxbNuM} 
\sum_{\alpha, \beta, \mu, \nu, j, k = 1}^{s}
\langle p_i, x_{\alpha} \rangle 
\Big( \sum_{N} \frac{x^N}{N!} b_{N + 1_{\alpha} + 1_{\beta}} \Big) 
\langle x_{\beta}, p_{j} \rangle 
\times \\ \times  
a_{1_{j} + 1_{k}} \langle p_{k}, x_{\mu} \rangle b_{1_{\mu} + 1_{\nu}} 
\Big( \sum_{M} \frac{x^{M}}{M!} u_{M} (\nu, \omega)  \Big) = 
\langle p_{i}, x_{\omega} \rangle,   
\end{multline}
for $i, \omega \in [s]$. 
Taking the leading term in $x$ yields: 
\begin{equation*}
\langle p_{i}, x_{\omega} \rangle = 
\sum_{\alpha, \beta, \mu, \nu, k, l = 1}^{s}
\langle p_{i}, x_{\alpha} \rangle b_{1_{\alpha} + 1_{\beta}} 
\langle x_{\beta}, p_{k} \rangle a_{1_{k} + 1_{l}} 
\langle p_{l}, x_{\mu} \rangle b_{1_{\mu} + 1_{\nu}} 
u_{\bar 0} (\nu, \omega), 
\end{equation*}
which is just 
$\langle p_{i}, x_{\omega} \rangle = 
\sum_{\mu, \nu = 1}^{s}
\langle p_{i}, x_{\mu} \rangle b_{1_{\mu} + 1_{\nu}} u_{\bar 0} (\nu, \omega)$, 
and one can see, that 
\begin{equation*} 
u_{\bar 0} (\alpha, \beta) = 
\sum_{k, l = 1}^{s}
\langle x_{\alpha}, p_{k} \rangle a_{1_{k} + 1_{l}} 
\langle p_{l}, x_{\beta} \rangle, 
\end{equation*}
for $\alpha, \beta \in [s]$. 
The equality 
(\ref{eq:pxbNuM}) 
linking $u_M (\alpha, \beta)$ and $b_{N}$ can be expressed as 
\begin{equation} 
\label{eq:baUemptyset}
\sum_{j, k = 1}^{s} (-1) B_{i, j} (x) a_{1_{j} + 1_{k}} U_{\emptyset} (x; k, \omega) = 0, 
\end{equation}
where $U_{\emptyset} (x; k, \omega) := 
\sum_{\alpha, \beta = 1}^{s}
\langle p_{k}, x_{\alpha} \rangle b_{1_{\alpha} + 1_{\beta}}
\sum_{N} (x^N/ N!) u_N (\beta, \omega)$, $k, \omega \in [s]$. 
Set 
$U_{k_1, \dots, k_r} (x; j, \omega) := 
\langle p_{k_{r}}, \dots \langle p_{k_{2}}, \langle p_{k_{1}}, U_{\emptyset} (x; j, \omega) 
\rangle \rangle \rangle$, and put 
\begin{equation*}
\begin{gathered}
U_{\emptyset}^{(0)} (j, \omega) := U_{\emptyset} (x; j, \omega)|_{x \to \bar 0}, \quad 
U_{k_{1}, \dots, k_{r}}^{(0)} (j, \omega) := 
U_{k_1, \dots, k_r} (x; j, \omega)|_{x \to \bar 0}. 
\end{gathered}
\end{equation*}
where $r \in \mathbb{Z}_{> 0}$, and $j, \omega, k_1, k_2, \dots, k_r \in [s]$. 
Recall, that it is assumed in this section that the braiding coefficients $q_{i', j'} = 1$, 
for all $i', j' \in [2 s]$. Simplifying the expressions, we obtain:  
\begin{equation} 
\label{eq:UemptysetUk}
\begin{gathered}
U_{\emptyset}^{(0)} (j, \omega) = \langle p_{j}, x_{\omega} \rangle, \\
U_{k_1, \dots, k_r}^{(0)} (j, \omega) = 
\sum_{\mu_1, \dots, \mu_r = 1}^{s} 
\Big( \prod_{i = 1}^{r} 
\langle p_{k_{i}}, x_{\mu_{i}} \rangle
\Big) 
\sum_{\alpha, \beta = 1}^{s}
\langle p_{j}, x_{\alpha} \rangle b_{1_{\alpha} + 1_{\beta}}
u_{1_{\mu_{1}} + \dots + 1_{\mu_{r}}} (\beta, \omega)
\end{gathered}
\end{equation}
Apply now the product $\prod_{i = 1}^{r} \langle p_{k_{i}}, - \rangle$ to the 
equation 
(\ref{eq:baUemptyset}) 
linking $B_{i, j} (x)$ and $U_{\emptyset} (x; l, \omega)$: 
\begin{equation*} 
\sum_{j, l \in [s]}
\sum_{I \subset [r]} 
\Big[ \Big( 
\prod_{i' \in I} \langle p_{k_{i'}}, - \rangle
\Big) B_{i, j} (x) \Big] 
a_{1_{j} + 1_{l}} 
\Big[ 
\Big( 
\prod_{i'' \in [r] \backslash I} 
\langle p_{k_{i''}}, - \rangle
\Big) U_{\emptyset} (x; l, \omega)
\Big] = 0. 
\end{equation*}
Evaluate this expression at $x \to \bar 0$.  
Since $B_{i, j}^{(0)} = - \sum_{\alpha, \beta \in [s]} 
\langle p_i, x_{\alpha} \rangle b_{1_{\alpha} + 1_{\beta}} \langle x_{\beta}, p_{j} \rangle$, 
the term 
$\sum_{i, j \in [s]} B_{i, j}^{(0)} a_{1_{j} + 1_{l}} U_{k_1, \dots, k_r}^{(0)} (l, \omega)$
corresponding to $I = \emptyset$ reduces just to 
$\sum_{\alpha, \beta, j, l \in [s]} (-1) \langle p_{i}, x_{\alpha} \rangle b_{1_{\alpha} + 1_{\beta}} 
\langle x_{\beta}, p_{j} \rangle a_{1_{j} + 1_{l}} 
U_{k_1, \dots, k_r}^{(0)} (l, \omega) = 
- U_{k_1, \dots, k_r}^{(0)} (l, \omega)$. 
Hence 
\begin{equation} 
\label{eq:U0B0_recursive}
U_{k_1, \dots, k_r}^{(0)} (i, \omega) = 
\sum_{j, l \in [s]}
\sum_{\substack{I \subset [r],\\ 
I \not = \emptyset}} B_{i, j, k_{I}}^{(0)} a_{1_{j} + 1_{l}} U_{k_{\ov{I}}}^{(0)} (l, \omega)
\end{equation}
for every $r \in \mathbb{Z}_{> 0}$ and $k_1, k_2, \dots, k_r \in [s]$, 
where one writes $B_{i, j, k_{I}}^{(0)}$ for $B_{i, j, l_1, \dots, l_{|I|}}^{(0)}$, if 
$I = \lbrace l_1 <l_2 < \dots < l_{|I|} \rbrace$, $|I|$ denotes the cardinality of $I$, 
$\ov{I} := [r] \backslash I$, and 
$U_{k_{\ov{I}}}^{(0)} (l, \omega)$ stands for 
$U_{l_{1}', \dots, l_{|\ov{I}|}'}^{(0)} (l, \omega)$, 
$\ov{I} = \lbrace l_{1}' < l_{2}' < \dots < l_{|\ov{I}|}' \rbrace$. 
Note, that if $I \not = \emptyset$, then $|\ov{I}| < r$. 
In terms of coefficients $\lbrace b_{N} \rbrace_{N}$, we have 
\begin{equation*} 
\begin{gathered} 
B_{i, j}^{(0)} = - \sum_{\alpha, \beta \in [s]} 
\langle p_i, x_{\alpha} \rangle b_{1_{\alpha} + 1_{\beta}} \langle x_{\beta}, p_{j} \rangle, \\
B_{i, j, k_1, \dots, k_r}^{(0)} = 
\sum_{\alpha, \beta, \mu_{1}, \dots, \mu_{r} \in [s]} 
\Big( 
\prod_{i' = 1}^{r} \langle p_{k_{i'}}, x_{\mu_{i'}} \rangle 
\Big) \langle p_{j}, x_{\beta} \rangle 
\langle p_{i}, x_{\alpha} \rangle 
b_{1_{\alpha} + 1_{\beta} + 1_{\mu_{1}} + \dots + 1_{\mu_{r}}}, 
\end{gathered}
\end{equation*}
for every $r \in \mathbb{Z}_{> 0}$ and $i, j, k_1, k_2, \dots, k_r \in [s]$. 
One can see that it is possible to express recursively 
all the quantities $U_{k_1, \dots, k_{r}}^{(0)} (l, \omega)$ in terms of 
$B_{i, j, k_1, \dots, k_r}^{(0)}$, $B_{i, j}^{(0)}$ and $a_{1_{i} + 1_{j}}$. 
Now the aim is to go to the quantities $A_{\alpha, \beta, \lambda_1, \dots, \lambda_r}^{(0)}$, 
$A_{\alpha, \beta}^{(0)}$. For them we have 
\begin{equation*} 
\begin{gathered}
A_{\alpha, \beta}^{(0)} = - \sum_{i, j \in [s]} \langle x_{\alpha}, p_{i} \rangle a_{1_{i} + 1_{j}} 
\langle p_{j}, x_{\beta}\rangle, \\
A_{\alpha, \beta, \lambda_1, \dots, \lambda_r}^{(0)} = 
\sum_{i, j, k_1, \dots, k_r \in [s]} 
\Big( 
\prod_{i' = 1}^{r} \langle x_{\lambda_{i'}}, p_{k_{i'}} \rangle
\Big) 
\langle x_{\beta}, p_{j} \rangle 
\langle x_{\alpha}, p_{i} \rangle  
a_{1_{i} + 1_{j} + 1_{k_{1}} + \dots + 1_{k_{r}}}, 
\end{gathered}
\end{equation*}
for every $r \in \mathbb{Z}_{> 0}$ and $\alpha, \beta, \lambda_1, \dots, \lambda_r \in [s]$. 
The link relating them to $U_{k_1, \dots, k_r}^{(0)} (l, \omega)$ is defined by 
(\ref{eq:a_uN}) 
and can be expressed as 
\begin{equation*} 
A_{\alpha, \beta} (p) = 
- \sum_{N} \frac{x^N}{N!} u_{N} (\alpha, \beta) \Big|_{x \to \langle x, \wt{S} (p) \rangle} = 
- \sum_{i, j \in [s]}
\langle x_{\alpha}, p_i \rangle a_{1_i + 1_j} 
U_{\emptyset} (x; j, \beta) |_{x \to \langle x, \wt{S} (p) \rangle}. 
\end{equation*}
One needs to take the multiple brackets $\prod_{i' = 1}^{r} \langle x_{\lambda_{i'}}, - \rangle$, 
and then to evaluate the left and the right-hand sides at $p \to \bar 0$. 
The left-hand side will yield just $A_{\alpha, \beta, \lambda_1, \dots, \lambda_{r}}^{(0)}$. 
It is convenient to introduce 
\begin{equation} 
\label{eq:Vemptyset}
V_{\emptyset} (p; \alpha, \beta) := 
\sum_{i, j \in [s]}
\langle x_{\alpha}, p_{i} \rangle a_{1_{i} + 1_{j}} 
U_{\emptyset} (x; j, \beta)|_{x \to \langle x, \wt{S} (p) \rangle}, 
\end{equation} 
and 
\begin{multline}
\label{eq:Vnu} 
V_{\nu_1, \dots, \nu_r} (p; \alpha, \beta) := 
\sum_{k_1, \dots, k_r \in [s]} \Big( 
\prod_{i' = 1}^{r} 
\sum_{\mu, i \in [s]}
b_{1_{\nu_{i'}} + 1_{\mu}} \langle x_{\mu}, p_{i} \rangle 
a_{1_{i} + 1_{k_{i'}}} \Big) 
\times \\ \times  
\sum_{j, l \in [s]}
\langle x_{\alpha}, p_{j} \rangle a_{1_{j} + 1_{l}} 
U_{k_1, \dots, k_r} (x; l, \beta)|_{x \to \langle x, \wt{S} (p) \rangle}, 
\end{multline}
for every $r \in \mathbb{Z}_{> 0}$ and $\alpha, \beta, \nu_1, \dots, \nu_r \in [s]$. 
In this notation we have 
\begin{equation*}
A_{\alpha, \beta} (p) = - V_{\emptyset} (p; \alpha, \beta), 
\end{equation*} 
so one needs to know what happens to $V_{\emptyset} (p; \alpha, \beta)$ once one 
applies $\prod_{i' = 1}^{r} \langle x_{\lambda_{i'}}, - \rangle$.  

\begin{prop} 
In the notation defined in 
(\ref{eq:Vemptyset}) 
and 
(\ref{eq:Vnu}), 
the following formulae are valid: 
\begin{equation*} 
\begin{gathered} 
\langle x_{\lambda}, V_{\emptyset} (p; \alpha, \beta) \rangle = 
- \sum_{\mu \in [s]} V_{\emptyset} (p; \lambda, \mu) V_{\mu} (p; \alpha, \beta), \\
\langle x_{\lambda}, V_{\nu_1, \dots, \nu_r} (p; \alpha, \beta) \rangle = 
- \sum_{\mu \in [s]} V_{\emptyset} (p; \lambda, \mu) V_{\nu_1, \dots, \nu_r, \mu} (p; \alpha, \beta), 
\end{gathered}
\end{equation*}
for $r \in \mathbb{Z}_{> 0}$, $\lambda, \alpha, \beta, \nu_1, \dots, \nu_r \in [s]$. 
\end{prop}
\noindent
\emph{Proof.} 
The computation is rather straightforward, but the formulae 
contain a lot of indices of summation. 
To save some letters of the Greek and Latin alphabets, let us 
simply write 
$\langle x_{\alpha}, p \rangle a_{1 + 1} \langle p, x_{\beta} \rangle$ 
for $\sum_{i, j \in [s]} \langle x_{\alpha}, p_{i} \rangle 
a_{1_{i} + 1_{j}} \langle p_{j}, x_{\beta} \rangle$, $\alpha, \beta \in [s]$. 
Similarly, 
$\langle p_{i}, x \rangle b_{1 + 1} \langle x, p_{j} \rangle$ will mean 
$\sum_{\alpha, \beta \in [s]} 
\langle p_{i}, x_{\alpha} \rangle b_{1_{\alpha} + 1_{\beta}} \langle x_{\beta}, p_{j} \rangle$, 
$i, j \in [s]$, etc. 
Expanding the definitions, we obtain: 
\begin{equation*} 
\begin{gathered}
V_{\emptyset} (p; \lambda, \mu) = 
\sum_{\sigma \in [s]} 
\langle x_{\lambda}, p \rangle a_{1 + 1} \langle p, x \rangle b_{1 + 1_{\sigma}} 
\sum_{N} \frac{\langle x, \wt{S} (p) \rangle^{N}}{N!} u_{N} (\sigma, \mu), \\
V_{\mu} (p; \alpha, \beta) = 
\sum_{k, l \in [s]} 
\big( b_{1_{\mu} + 1} \langle x, p \rangle a_{1 + 1_{k}} \big)   
\big[ \langle x_{\alpha}, p \rangle a_{1 + 1_{l}} \big]
U_{k} (x; l, \beta)|_{x \to \langle x, \wt{S} (p) \rangle}, 
\end{gathered}
\end{equation*}
where $\lambda, \mu, l, \beta \in [s]$, and for $U_{k} (x; l, \beta)$, $k \in [s]$, 
we have 
\begin{equation*} 
U_{k} (x; l, \beta) = 
\sum_{\sigma, \rho \in [s]} 
\langle p_{k}, x_{\rho} \rangle 
\big( \langle p_{l}, x \rangle b_{1 + 1_{\sigma}} \big) 
\sum_{N} \frac{x^N}{N!} u_{N + 1_{\rho}} (\sigma, \beta). 
\end{equation*} 
This assembles as follows into a formula for $V_{\mu} (p; \alpha, \beta)$: 
\begin{multline*} 
V_{\mu} (p; \alpha, \beta) = 
\sum_{\sigma, \rho \in [s]} 
\big( b_{1_{\mu} + 1} \langle x, p \rangle a_{1 + 1} \langle p, x_{\rho} \rangle \big)
\times \\ \times
\big[ \langle x_{\alpha}, p \rangle a_{1 + 1} \langle p, x \rangle b_{1 + 1_{\sigma}} \big] 
\sum_{N} 
\frac{\langle x, \wt{S} (p)\rangle^{N}}{N!} u_{N + 1_{\rho}} (\sigma, \beta). 
\end{multline*}
Computing the bracket $\langle x_{\lambda}, - \rangle$ on $V_{\emptyset} (p; \alpha, \beta)$ yields: 
\begin{multline*} 
\langle x_{\lambda}, V_{\emptyset} (p; \alpha, \beta) \rangle = 
- \sum_{\rho, \pi \in [s]} 
\big[ 
\langle x_{\alpha}, p \rangle a_{1 + 1} \langle p, x \rangle b_{1 + 1_{\rho}} 
\big] 
\times \\ \times
\sum_{N} \frac{\langle x, \wt{S} (p) \rangle^{N}}{N!} 
u_{N + 1_{\pi}} (\rho, \beta) 
\langle \langle x_{\pi}, \wt{S} (p) \rangle, x_{\lambda} \rangle. 
\end{multline*}
The double bracket in the right-hand side is of the shape:  
\begin{multline*} 
\langle \langle x_{\pi}, \wt{S} (p) \rangle, x_{\lambda} \rangle = 
- A_{\lambda, \pi} (p) = 
\sum_{N} \frac{\langle x, \wt{S} (p) \rangle^{N}}{N!} u_{N} (\lambda, \pi) 
= \\ =  
- \sum_{\sigma, \mu \in [s]} 
\big( 
\langle x_{\lambda}, p \rangle a_{1 + 1} \langle p, x \rangle b_{1 + 1_{\sigma}} 
\big) 
\Big[ 
\sum_{N} \frac{\langle x, \wt{S} (p) \rangle^{N}}{N!} u_{N} (\sigma, \mu) 
\Big] 
\big(
b_{1_{\mu} + 1} \langle x, p \rangle a_{1 + 1} \langle p, x_{\pi} \rangle
\big). 
\end{multline*}
Substituting this expression into the previous formula, we obtain just 
\begin{equation*} 
\langle x_{\lambda}, V_{\emptyset} (p; \alpha, \beta) \rangle = 
- \sum_{\mu \in [s]} V_{\emptyset} (p; \lambda, \mu) V_{\mu} (p; \alpha, \beta),  
\end{equation*}
as claimed. 

The computations for $V_{\nu_1, \dots, \nu_r} (p; \alpha, \beta)$ are totally similar.  
Expanding the definitions, one obtains: 
\begin{multline*} 
V_{\nu_1, \dots, \nu_r} (p; \alpha, \beta) = 
\sum_{\sigma, \rho_1, \dots, \rho_r \in [s]} 
\Big\lbrace
\prod_{i' = 1}^{r}
\big( b_{1_{\nu_{i'}} + 1} \langle x, p \rangle a_{1 + 1} \langle p, x_{\rho_{i'}} \rangle \big)
\Big\rbrace 
\times \\ \times
\big[ \langle x_{\alpha}, p \rangle a_{1 + 1} \langle p, x \rangle b_{1 + 1_{\sigma}} \big] 
\sum_{N} 
\frac{\langle x, \wt{S} (p)\rangle^{N}}{N!} u_{N + 1_{\rho_{1}} + \dots + 1_{\rho_{r}}} (\sigma, \beta),  
\end{multline*} 
where $\alpha, \beta, \nu_1, \dots, \nu_r \in [s]$. 
Taking the bracket $\langle x_{\lambda}, - \rangle$ causes an additional shift 
by $1_{*}$ in the lower index of $u$ in the right-hand side, resulting in  
\begin{multline*} 
\langle x_{\lambda}, V_{\nu_1, \dots, \nu_r} (p; \alpha, \beta) \rangle = 
- \sum_{\pi, \sigma, \rho_1, \dots, \rho_r \in [s]} 
\Big\lbrace
\prod_{i' = 1}^{r}
\big( b_{1_{\nu_{i'}} + 1} \langle x, p \rangle a_{1 + 1} \langle p, x_{\rho_{i'}} \rangle \big)
\Big\rbrace 
\times \\ \times
\big[ \langle x_{\alpha}, p \rangle a_{1 + 1} \langle p, x \rangle b_{1 + 1_{\sigma}} \big] 
\sum_{N} 
\frac{\langle x, \wt{S} (p)\rangle^{N}}{N!} 
u_{N + 1_{\rho_{1}} + \dots + 1_{\rho_{r}} + 1_{\pi}} (\sigma, \beta)
\langle \langle x_{\pi}, \wt{S} (p) \rangle, x_{\lambda} \rangle. 
\end{multline*} 
It remains to substitute the expression for the 
double bracket $\langle \langle x_{\pi}, \wt{S} (p) \rangle, x_{\lambda} \rangle$ 
and to simplify the result. 
This yields the second formula claimed. 
\qed 

Now we can return to the equation 
$A_{\alpha, \beta} (p) = - V_{\emptyset} (p; \alpha, \beta)$. 
From the just proved proposition, it follows, that 
the brackets 
$\langle p_{k_r}, \dots \langle p_{k_2}, \langle p_{k_1}, A_{\alpha, \beta} (p) 
\rangle \rangle \dots \rangle$, 
$r \in \mathbb{Z}_{>0}$, $\alpha, \beta, k_1, \dots, k_r \in [s]$, 
can be expressed in terms of 
$V_{\emptyset} (p; \alpha, \beta)$ and 
the higher order coefficients 
$V_{k_1, \dots, k_r} (p; \alpha, \beta)$. 
For example, 
\begin{multline*} 
\langle x_{\lambda_2}, \langle x_{\lambda_1}, - A_{\alpha, \beta} (p) \rangle \rangle = 
- \sum_{\mu_{1} \in [s]} 
\Big\lbrace 
\langle x_{\lambda_2}, V_{\emptyset} (p; \lambda_1, \mu_1) \rangle 
V_{\mu_1} (p; \alpha, \beta) 
+ \\ + 
V_{\emptyset} (p; \lambda_1, \mu_1) 
\langle 
x_{\lambda_2}, V_{\mu_1} (p; \alpha, \beta)
\rangle
\Big\rbrace = 
(- 1)^{2}
\sum_{\mu_1, \mu_2 \in [s]} 
\Big\lbrace 
V_{\emptyset} (p; \lambda_2, \mu_2) V_{\mu_2} (p; \lambda_1, \mu_1) 
\times \\ \times 
V_{\mu_1} (p; \alpha, \beta) +  
V_{\emptyset} (p; \lambda_2, \mu_2) V_{\emptyset} (p; \lambda_1, \mu_1) V_{\mu_1, \mu_2} (p; \alpha, \beta)
\Big\rbrace. 
\end{multline*} 
One can obtain a generic formula for any $r$ by induction. 
Extracting the coefficient corresponding to $p^{\bar 0}$, 
and introducing the notation 
$V_{\emptyset}^{(0)} (\alpha, \beta) := V_{\emptyset} (p; \alpha, \beta)|_{p \to \bar 0}$ and 
$V_{\nu_1, \dots, \nu_{r}}^{(0)} (\alpha, \beta) := 
V_{\nu_1, \dots, \nu_{r}} (p; \alpha, \beta)|_{p \to \bar 0}$, 
$\alpha, \beta, \nu_1, \dots, \nu_r \in [s]$, $r \in \mathbb{Z}_{>0}$, 
one obtains: 
\begin{equation} 
\label{eq:AzeroVzero}
A_{\lambda_0, \mu_0, \lambda_1, \dots, \lambda_r}^{(0)} = 
(- 1)^{r - 1} 
\sum_{\substack{f : \lbrace 0, \dots, r \rbrace \to 
\lbrace 0, \dots, r - 1 \rbrace, \\
0 \leqslant f (j') \leqslant j' - 1, \\
j' = 1, \dots, r 
}}
\sum_{\mu_1, \dots, \mu_r \in [s]} 
\prod_{i' = 0}^{r} 
V_{\mu_{f^{- 1} (i')}}^{(0)} (\lambda_{i'}, \mu_{i'}), 
\end{equation}
for any $r \in \mathbb{Z}_{>0}$, any $\lambda_0, \mu_0, \lambda_1, \dots, \lambda_r \in [s]$, 
and  
$V_{\mu_{f^{- 1} (i')}}^{(0)} (\lambda_{i'}, \mu_{i'})$ is understood as 
$V_{\mu_{t_{1} [i', f]}, \dots, \mu_{t_{m} [i', f]}}^{(0)} (\lambda_{i'}, \mu_{i'})$, 
where $\lbrace t_1 [i', f] < \dots < t_m [i', f] \rbrace$ is the 
coimage of an element $i' \in \lbrace 0, 1, \dots, r \rbrace$ under the map $f$. 
In terms of the coefficients $U_{\emptyset}^{(0)} (j, \omega) = \langle p_j, x_{\omega} \rangle$, 
and $U_{k_1, \dots, k_r}^{(0)} (j, \omega)$, $r \geqslant 1$, one obtains: 
\begin{equation*}
V_{\emptyset}^{(0)} (\alpha, \beta) = 
\sum_{i, j \in [s]}
\langle x_{\alpha}, p_{i} \rangle a_{1_i + 1_j} 
\langle p_{j}, x_{\beta} \rangle, 
\end{equation*}
and 
\begin{multline*} 
V_{\nu_1, \dots, \nu_r}^{(0)} (\alpha, \beta) = 
\sum_{k_1, \dots, k_r \in [s]} \Big( 
\prod_{i' = 1}^{r} 
\sum_{\mu, i \in [s]}
b_{1_{\nu_{i'}} + 1_{\mu}} \langle x_{\mu}, p_{i} \rangle 
a_{1_{i} + 1_{k_{i'}}} \Big) 
\times \\ \times  
\sum_{j, l \in [s]}
\langle x_{\alpha}, p_{j} \rangle a_{1_{j} + 1_{l}} 
U_{k_1, \dots, k_r}^{(0)} (l, \beta),   
\end{multline*}
for $r \in \mathbb{Z}_{>0}$ and $\alpha, \beta, \nu_1, \dots, \nu_r \in [s]$. 

In other words, we can link the coefficients 
$A_{\alpha, \beta, \lambda_1, \dots, \lambda_r}^{(0)}$ and 
$B_{i, j, k_1, \dots, k_{r'}}^{(0)}$ ($r, r' \in \mathbb{Z}_{> 0}$), in three steps: 
1) express recursively $U_{*}^{(0)}$ via $B_{*}^{(0)}$ 
(see (\ref{eq:UemptysetUk}) and (\ref{eq:U0B0_recursive})); 
2) express $V_{*}^{(0)}$ via $U_{*}^{(0)}$; 
3) express $A_{*}^{(0)}$ via $V_{*}^{(0)}$. 
The resulting formulae describing the transition $B_{*}^{(0)} \to A_{*}^{(0)}$ 
is a noncommutative generalization of the 
Legendre transformation in case all braidings $q_{i', j'} = 1$, $i', j' = 1, 2, \dots, 2 s$.

\section{$q$-commutative units of measurement.} 
In this section we describe the $q$-Legendre transformation in a generic case 
by introducing the \emph{``$q$-commutative units of measurement''}. 
The symbol $q$ stands for a $2 s \times 2 s$ matrix of formal variables $q_{i', j'}$ 
satisfying $q_{i', i'} = 1$, and $q_{i', j'} q_{j', i'} = 1$, $i', j' = 1, 2, \dots, 2 s$. 

In physics, one can not add quantities which have different units of measurements. 
One can not add 5 grams with 10 centimetres. 
At the same time, if we have a mass $M$ of something, and a length $L$ of something else, 
then fixing the units of measurement, say $M_0 = 1 \mathit{gram}$ and $L_{0} = 1 \mathit{cm}$, 
gives an opportunity to consider a sum $M/ M_{0} + L/ L_{0}$. 
Having in mind this analogy, let us describe the general set-up. 
We have a collection $p_1, p_2, \dots, p_s, x_1, x_2, \dots, x_s$ of 
$q$-commuting variables 
\begin{equation*} 
z_{i'} z_{j'} = q_{j', i'} z_{j'} z_{i'}, \quad 
z_{k} := p_{k}, \quad z_{s + \alpha} := x_{\alpha},  
\end{equation*} 
where $i', j' \in [2 s] = \lbrace 1, 2, \dots, 2 s \rbrace$, 
$k, \alpha \in [s] = \lbrace 1, 2, \dots, s \rbrace$. 
Denote $\mathcal{A}_{q}$ the quantum affine space generated by $z_{1}, z_{2}, \dots, z_{2 s}$. 

Next, adjoin to $\mathcal{A}_{q}$ a collection of variables 
$h_{\alpha, i}$, $\alpha, i \in [s]$, which ``behave like'' products $x_{\alpha} p_{i}$, i.e. 
\begin{equation*} 
h_{\alpha, i} z_{k'} = 
q_{k', s + \alpha} q_{k', i} 
z_{k'} h_{\alpha, i}, \quad 
h_{\alpha, i} h_{\beta, j} = 
q_{s + \beta, s + \alpha} q_{s + \beta, i} q_{j, s + \alpha} q_{j, i} 
h_{\beta, j} h_{\alpha, i}, 
\end{equation*}
where $\alpha, \beta, i, j \in [s]$, and $k' \in [2 s]$. 
Denote the resulting quantum affine space $\wt{\mathcal{A}}_{q}$. 
For a multi-index $N \in \mathbb{Z}_{\geqslant 0}^{s}$ of length $s$, 
$N = (N_1, N_2, \dots, N_s)$, let us write 
\begin{equation*} 
p^{N} := p_{s}^{N_s} \dots p_{1}^{N_1}, \quad 
x^{N} := x_{s}^{N_s} \dots x_{1}^{N_1},
\end{equation*}
It is convenient to have a notation $Q_{(K, L), (M, N)}$, 
where $K, L, M, N \in \mathbb{Z}_{\geqslant 0}^{s}$, for the coefficient 
defined by 
\begin{equation*} 
x^{M} p^{N} x^{K} p^{L} = Q_{(K, L), (M, N)} x^{K} p^{L} x^{M} p^{N}. 
\end{equation*}
Then when we say, that a quantity $f$ ``behaves like'' $x^M p^N$, while 
a quantity $g$ ``behaves like'' $x^K p^L$, has just a meaning that 
$f g = Q_{(K, L), (M, N)} g f$. Explicitly: 
\begin{equation*} 
Q_{(K, L), (M, N)} := 
\Big( 
\prod_{\beta, i = }^{s} q_{s + \beta, i}^{K_{\beta} N_{i}}
\Big)
\Big[ 
\prod_{\beta, \alpha = 1}^{s} q_{s + \beta, s + \alpha}^{K_{\beta} M_{\alpha}}
\Big] 
\Big( 
\prod_{j, i = 1}^{s} q_{j, i}^{L_j N_i}
\Big) 
\Big[ 
\prod_{j, \alpha = 1}^{s} q_{j, s + \alpha}^{- L_{j} M_{\alpha}}
\Big]. 
\end{equation*}
This naturally extends to any multi-indices $K, L, M, N \in \mathbb{Z}^s$, 
not necessary in $\mathbb{Z}_{\geqslant 0}^{s}$, 
so it makes sense to speak $f$ ``behaves like'' $x^{M} p^{N}$ for any $M, N \in \mathbb{Z}^s$. 
The quantum affine space $\wt{\mathcal{A}}_{q}$ generated by 
$p_{i}$, $x_{\alpha}$, and $h_{\alpha, i}$, $i, \alpha \in [s]$, 
is equipped with a bracket $\langle -, - \rangle: 
\wt{\mathcal{A}}_{q} \times \wt{\mathcal{A}}_{q} \to 
\wt{\mathcal{A}}_{q}$ defined as a bilinear map by  
\begin{equation*} 
\langle p_{i}, x_{\alpha} \rangle = h_{\alpha, i}, \quad 
\langle x_{\alpha}, p_{i} \rangle = - q_{i, s + \alpha} h_{\alpha, i}, \quad 
\langle p_i, p_j \rangle = 0, \quad 
\langle x_{\alpha}, x_{\beta} \rangle = 0, 
\end{equation*} 
and 
\begin{multline*}
\langle z_{k_r} \dots z_{k_1}, z_{l_{r'}} \dots z_{l_1} \rangle = 
\sum_{i' = 1}^{r} \sum_{j' = 1}^{r'} 
z_{k_r} \dots \check z_{k_{i'}} \dots z_{k_1} 
\langle z_{k_{i'}}, z_{l_{j'}} \rangle 
\times \\ \times 
z_{l_{r'}} \dots \check z_{l_{j'}} \dots z_{l_1}
\Big( \prod_{i'' = 1}^{i' - 1} q_{i'', i'} \Big) 
\prod_{j'' = j' + 1}^{r'} 
 q_{j', j''}, 
\end{multline*}
where $i, j, \alpha, \beta \in [s]$, 
$r, r' \in \mathbb{Z}_{> 0}$, $k_1, \dots, k_r, l_1, \dots, l_{r'} \in [2 s]$. 

The next step is to adjoin to $\wt{A}_{q}$ collections of variables 
$\bar{a}_{i, j}$ and $\bar b_{\alpha, \beta}$, $i, j, \alpha, \beta \in [s]$, 
which behave like $p^{-(1_{i} + 1_{j})}$ and $x^{- (1_{\alpha} + 1_{\beta})}$, respectively. 
Write symbolically 
\begin{equation*} 
\bar a_{i, j} \Longrightarrow p^{-(1_{i} + 1_{j})}, \quad 
\bar b_{\alpha, \beta} \Longrightarrow x^{- (1_{\alpha} + 1_{\beta})}, 
\end{equation*}
to express this fact. 
In this notation, we already have 
\begin{equation*} 
p_{i} \Longrightarrow p^{1_{i}}, \quad 
x_{\alpha} \Longrightarrow x^{1_{\alpha}}, \quad 
h_{\alpha, i} \Longrightarrow x^{1_{\alpha}} p^{1_{i}}, 
\end{equation*}
for $i, \alpha \in [s]$. 
So we obtain yet another quantum affine space, denote it $\wh{\mathcal{A}}_{q}$, 
with a finite number of generators 
$p_i$, $x_{\alpha}$, $h_{\alpha, i}$, $\bar a_{i, j}$, $\bar b_{\alpha, \beta}$, 
where $i, j, \alpha, \beta = 1, 2, \dots, s$. 

Now we need to construct from $\wh{\mathcal{A}}_{q}$ two other quantum affine spaces, 
$\wh{\mathcal{A}}_{q}^{(\mathit{coor})}$ and $\wh{\mathcal{A}}_{q}^{(\mathit{mom})}$, 
this time having \emph{infinite} numbers of generators each.  
The first one $\wh{\mathcal{A}}_{q}^{(\mathit{coor})}$ is obtained 
by adjoining 
\begin{equation*} 
\bar B_{i, j, k_1, k_2, \dots, k_{r}} 
\Longrightarrow p^{1_{i} + 1_{j} + 1_{k_1} + 1_{k_2} + \dots + 1_{k_r}}, 
\end{equation*}
where $r \in \mathbb{Z}_{> 0}$, and $i, j, k_1, k_2, \dots k_r \in [s]$. 
The second one $\wh{\mathcal{A}}_{q}^{(\mathit{mom})}$ is totally similar, and is obtained from 
$\wh{\mathcal{A}}_{q}$ by adjoining 
\begin{equation*} 
\bar A_{\alpha, \beta, \lambda_1, \lambda_2, \dots, \lambda_r} \Longrightarrow 
x^{1_{\alpha} + 1_{\beta} + 1_{\lambda_1} + 1_{\lambda_2} + \dots + 1_{\lambda_r}}, 
\end{equation*}
where $r \in \mathbb{Z}_{> 0}$, and $\alpha, \beta, \lambda_1, \lambda_2, \dots, \lambda_r \in [s]$.

Factor out an ideal $\mathcal{I}_{q}$ in $\wh{\mathcal{A}}_{q}$ generated by relations 
\begin{equation}
\label{eq:idealIq_rels} 
\begin{gathered}
\sum_{\alpha, \beta, i, j \in [s]} 
\langle p_{l}, x_{\alpha} \rangle \bar b_{\alpha, \beta} 
\langle x_{\beta}, p_{i} \rangle \bar a_{i, j} 
\langle p_{j}, x_{\omega} \rangle = \langle p_{l}, x_{\omega} \rangle,\\ 
\sum_{\alpha, \beta, i, j \in [s]} 
\langle x_{\omega}, p_{i} \rangle
\bar a_{i, j}
\langle p_{j}, x_{\alpha} \rangle 
\bar b_{\alpha, \beta} 
\langle x_{\beta}, p_{l} \rangle = 
\langle x_{\omega}, p_{l} \rangle,\\ 
\bar a_{i, j} = Q_{(\bar 0, - 1_{j}), (\bar 0, - 1_{i})} \bar a_{j, i}, \quad 
\bar b_{\alpha, \beta} = 
Q_{(- 1_{\beta}, \bar 0), (- 1_{\alpha}, \bar 0)}
\bar b_{\beta, \alpha}, 
\end{gathered}
\end{equation}
where $\alpha$, $\beta$, $\omega$, $i$, $j$, $l$ vary over $[s]= \lbrace 1, 2, \dots, s\rbrace$. 
Denote the result $\mathcal{B}_{q} := \mathcal{A}_{q}/ \mathcal{I}_{q}$, and 
keep the symbols 
$x_{\alpha}$, $p_{i}$, $h_{\alpha, i}$, $\bar a_{i, j}$, and $\bar b_{\alpha, \beta}$, 
to denote the equivalence classes 
$[x_{\alpha}]$, $[p_{i}]$, $[h_{\alpha, i}]$, $[\bar a_{i, j}]$, and $[\bar b_{\alpha, \beta}]$ 
in $\mathcal{B}_{q}$, respectively. 
Similarly, factor out an ideal $\mathcal{I}_{q}^{(\mathit{coor})}$ in 
$\wh{\mathcal{A}}_{q}^{(\mathit{coor})}$, 
$\mathcal{B}_{q}^{(\mathit{coor})} := \wh{\mathcal{A}}_{q}^{(\mathit{coor})}/ 
\mathcal{I}_{q}^{(\mathit{coor})}$, defined 
by these relations and the relations 
\begin{equation*} 
\bar B_{k_{\sigma (1)}, k_{\sigma (2)}, \dots, k_{\sigma(r + 2)}} = 
\Big( 
\prod_{\substack{1 \leqslant i' < j' \leqslant r + 2, \\
\sigma^{-1} (i') > \sigma^{-1} (j')
}} q_{k_{k'}, k_{i'}} 
\Big) \bar B_{k_{1}, k_{2}, \dots, k_{r + 2}}, 
\end{equation*} 
for any $\sigma \in S_{r + 2}$
(the symmetric group on $r + 2$ symbols), 
$r \in \mathbb{Z}_{> 0}$, 
and the indices $k_1, k_2, \dots, k_{r + 2} \in [s]$. 
In analogy with $\wh{\mathcal{A}}_{q}^{(\mathit{coor})}$, look at 
$\wh{\mathcal{A}}_{q}^{(\mathit{mom})}$, 
define an ideal $\mathcal{I}_{q}^{(\mathit{mom})}$ by 
relations 
(\ref{eq:idealIq_rels})
(as for the ideal $\mathcal{I}_{q}$) and the relations 
\begin{equation*} 
\bar A_{\lambda_{\sigma (1)}, \lambda_{\sigma (2)}, \dots, \lambda_{\sigma (r + 2)}} = 
\Big( 
\prod_{\substack{1 \leqslant i' < j' \leqslant r + 2,\\ 
\sigma^{-1} (i') > \sigma^{-1} (j')}}
q_{s + \lambda_{j'}, s + \lambda_{i'}}
\Big) 
\bar A_{\lambda_1, \lambda_2, \dots, \lambda_{r + 2}}, 
\end{equation*}
where $\sigma \in S_{r + 2}$, $r \in \mathbb{Z}_{> 0}$, 
$\lambda_1, \lambda_2, \dots, \lambda_{r + 2} \in [s]$. 
Set $\mathcal{B}_{q}^{(\mathit{mom})} := \wh{\mathcal{A}}_{q}^{(\mathit{mom})}/ 
\mathcal{I}_{q}^{(\mathit{mom})}$. 
We accept the convention not to write the square brackets for the 
canonical images of the generators of $\wh{\mathcal{A}}_{q}^{(\mathit{coor})}$ and 
$\wh{\mathcal{A}}_{q}^{(\mathit{mom})}$ in $\mathcal{B}_{q}^{(\mathit{coor})}$ and 
$\mathcal{B}_{q}^{(\mathit{mom})}$, 
respectively. 

Let $\mathcal{T}_{q}^{(\mathit{coor})} \subset \mathcal{B}_{q}^{(\mathit{coor})}$ 
be the the subalgebra of $\mathcal{B}_{q}^{(\mathit{coor})}$ 
formed by all linear combinations of monomials in $h_{\alpha, i}$, 
$\bar a_{i, j}$, $\bar b_{\alpha, \beta}$, and $\bar B_{i, j, k_1, \dots, k_r}$ 
(i.e. no generators $p_i$ or $x_{\alpha}$).  
Similarly, let $\mathcal{T}_{q}^{(\mathit{mom})} \subset \mathcal{B}_{q}^{(\mathit{mom})}$ be the 
subalgebra spanned over the monomials in 
$h_{\alpha, i}$, $\bar a_{i, j}$, $\bar b_{\alpha, \beta}$, and 
$\bar A_{\alpha, \beta, \lambda_1, \dots, \lambda_r}$ (no $p_i$ or $x_{\alpha}$). 
We suggest to define the \emph{$q$-Legendre transformation} as a map 
\begin{equation*} 
L_{q}: \mathcal{T}_{q}^{(\mathit{mom})} \to \mathcal{T}_{q}^{(\mathit{coor})}. 
\end{equation*}
Note, the the ``direction'' above is from ``momenta'' to ``coordinates'', 
which is dual to the direction at the 
classical (commutative) level:  
the map $L_q$ should be perceived as an analogue of going from 
$S (x)$ (classical action as a function of coordinates) to 
$\wt{S} (p)$ (classical action as a function of momenta). 
Essentially, we have already described the map $L_{q}$ 
corresponding to $q_{i', j'} = 1$, $i', j' \in [2 s]$, 
in the previous section. 
It is necessary to replace 
$a_{1_i + 1_j}$ with $\bar a_{i, j}$, 
$b_{1_{\alpha} + 1_{\beta}}$ with $\bar b_{\alpha, \beta}$, 
$A_{\alpha, \beta, \lambda_1, \dots, \lambda_r}^{(0)}$ with 
$\bar A_{\alpha, \beta, \lambda_1, \dots, \lambda_r}$, and 
$B_{i, j, k_{1}, \dots, k_{r}}^{(0)}$ with 
$\bar B_{i, j, k_{1}, \dots, k_{r}}$ in all the formulae. 
This defines a map 
$L_{\mathbf{1}}: \mathcal{T}_{\mathbf{1}}^{(\mathit{mom})} \to \mathcal{T}_{\mathbf{1}}^{(\mathit{coor})}$, 
where the index $\mathbf{1}$ denotes the $2 s \times 2 s$ matrix with the entries all equal to 1.  
Moreover, this map is an isomorphism of algebras, such that 
each generator $h_{\alpha, i}, \bar a_{i, j}, 
\bar A_{\alpha, \beta, \lambda_1, \dots, \lambda_r} \in \mathcal{T}_{q}^{(\mathit{mom})}$
is mapped to a polynomial in 
generators $h_{\alpha', i'}, \bar b_{\alpha', \beta'}, 
\bar B_{i', j', k_1', \dots, k_r'} \in \mathcal{T}_{q}^{(\mathit{coor})}$. 
\begin{thm} 
There exists a canonical extension $L_{q}: \mathcal{T}_{q}^{(\mathit{mom})} \to 
\mathcal{T}_{q}^{(\mathit{coor})}$ of 
the map $L_{\mathbf{1}}: \mathcal{T}_{\mathbf{1}}^{(\mathit{mom})} \to 
\mathcal{T}_{\mathbf{1}}^{(\mathit{coor})}$, 
analytical in $\lbrace q_{i', j'} \rbrace_{i', j' \in [2 s]}$, establishing an algebra isomorphism. 
\end{thm}
\noindent
\emph{Proof.} 
We describe the map $L_{q}$ explicitly. 
For that we will need to introduce what one terms the ``$q$-commutative units of measurement''. 
They allow to modify systematically the coefficients in the formulae for $L_{\mathbf{1}}$ 
in order to produce $L_{q}$. Recall, that we have introduced the algebras with the following generators: 
\begin{equation*} 
\begin{aligned} 
\mathcal{A}_{q}:\quad &x_{\alpha}, \quad p_{i},\\ 
\wt{\mathcal{A}}_{q}:\quad &x_{\alpha}, \quad p_{i}, \quad h_{\alpha, i},\\ 
\wh{\mathcal{A}}_{q}:\quad &x_{\alpha}, \quad p_{i}, \quad h_{\alpha, i}, \quad 
\bar a_{i, j}, \quad \bar b_{\alpha, \beta},\\  
\wh{\mathcal{A}}_{q}^{(\mathit{coor})}:\quad &x_{\alpha}, \quad p_{i}, \quad h_{\alpha, i}, \quad 
\bar a_{i, j}, \quad \bar b_{\alpha, \beta},\quad 
\bar B_{i, j, k_1, \dots, k_r},\\  
\wh{\mathcal{A}}_{q}^{(\mathit{mom})}:\quad &x_{\alpha}, \quad p_{i}, \quad h_{\alpha, i}, \quad 
\bar a_{i, j}, \quad \bar b_{\alpha, \beta},\quad 
\bar A_{\alpha, \beta, \lambda_1, \dots, \lambda_r},\\  
\mathcal{T}_{q}^{(\mathit{\mathit{coor}})}:\quad &h_{\alpha, i}, \quad 
\bar a_{i, j}, \quad \bar b_{\alpha, \beta},\quad 
\bar B_{i, j, k_1, \dots, k_r},\\  
\mathcal{T}_{q}^{(\mathit{\mathit{mom}})}:\quad &h_{\alpha, i}, \quad 
\bar a_{i, j}, \quad \bar b_{\alpha, \beta},\quad 
\bar A_{\alpha, \beta, \lambda_1, \dots, \lambda_r},
\end{aligned}
\end{equation*}
where $i, j, k_1, \dots, k_r, \alpha, \beta, \lambda_1, \dots, \lambda_r$ vary over $[s]$, and 
$r$ varies over $\mathbb{Z}_{> 0}$. 
One naturally extends the bracket $\langle -, - \rangle: \wt{\mathcal{A}}_{q} \times 
\wt{\mathcal{A}}_{q} \to \wt{\mathcal{A}}_{q}$ to the algebras 
$\wh{\mathcal{A}}_{q}$, $\wh{\mathcal{A}}_{q}^{(\mathit{coor})}$, and 
$\wh{\mathcal{A}}_{q}^{(\mathit{mom})}$, 
assuming the properties 
$\langle g F_1, F_2 \rangle = g \langle F_1, F_2 \rangle$ and 
$\langle F_1, F_2 g \rangle = \langle F_1, F_2 \rangle g$, whenever 
$g$ is one of the generators mentioned which is different from $p_{i}$ or $x_{\alpha}$
(where $(F_1, F_2)$ is an arbitrary pair of elements of one of the three algebras mentioned).  
Consider now additional symbols 
$\varkappa_1, \varkappa_2, \dots, \varkappa_{s}$ and 
$\theta_1, \theta_2, \dots, \theta_{s}$ with the properties 
\begin{equation*} 
\varkappa_{i} \Longrightarrow p^{- 1_{i}}, \quad 
\theta_{\alpha} \Longrightarrow x^{- 1_{\alpha}}, 
\end{equation*} 
for $i, \alpha \in [s]$. Adjoining them to the algebras 
$\wh{\mathcal{A}}_{q}^{(\mathit{coor})}$ and $\wh{\mathcal{A}}_{q}^{(\mathit{mom})}$, 
and then localizing with respect to $\varkappa_{i}$, $\theta_{\alpha}$ ($i, \alpha \in [s]$), 
we obtain two more algebras: 
\begin{equation*} 
\begin{aligned}
\wh{\mathcal{A}}_{q}^{(\mathit{coor, ext})}:\quad 
&x_{\alpha}, \quad p_{i}, \quad h_{\alpha, i}, \quad 
\bar a_{i, j}, \quad \bar b_{\alpha, \beta}, \quad  
\bar B_{i, j, k_1, \dots, k_r}, \quad 
\varkappa_{i}, \quad \theta_{\alpha},
\quad 
\varkappa_{i}^{-1}, \quad \theta_{\alpha}^{-1},
\\ 
\wh{\mathcal{A}}_{q}^{(\mathit{mom, ext})}:\quad &x_{\alpha}, \quad p_{i}, \quad h_{\alpha, i}, \quad 
\bar a_{i, j}, \quad \bar b_{\alpha, \beta},\quad 
\bar A_{\alpha, \beta, \lambda_1, \dots, \lambda_r}, \quad \varkappa_{i}, \quad \theta_{\alpha}
\quad 
\varkappa_{i}^{-1}, \quad \theta_{\alpha}^{-1}.  
\end{aligned}
\end{equation*}
Since in these algebras one has 
\begin{equation*} 
\varkappa_{i} p_{i} \Longrightarrow x^{\bar 0} p^{\bar 0}, \quad 
\theta_{\alpha} x_{\alpha} \Longrightarrow x^{\bar 0} p^{\bar 0}, 
\end{equation*}
where $i, \alpha \in [s]$, $\bar 0 = (0, 0, \dots 0)$ (length $s$), 
intuitively, it is natural to perceive the multiplication 
by $\varkappa_{i}$ as dividing by the ``unit of measurement'' of the momentum $p_i$, 
and the multiplication by $\theta_{\alpha}$ as dividing by the 
``unit of measurement'' of the coordinate $x_{\alpha}$. 
This time these units are not commutative, but $q$-commutative.  
There are canonical embeddings 
$\wh{\mathcal{A}}_{q}^{(\mathit{coor})} \subset 
\wh{\mathcal{A}}_{q}^{(\mathit{coor, ext})}$ and 
$\wh{\mathcal{A}}_{q}^{(\mathit{mom})} \subset 
\wh{\mathcal{A}}_{q}^{(\mathit{mom, ext})}$, 
and one can also find a copy of $\wh{\mathcal{A}}_{\mathbf{1}}^{(\mathit{coor})}$ inside 
$\wh{\mathcal{A}}_{q}^{(\mathit{coor, ext})}$, and a copy of 
$\wh{\mathcal{A}}_{\mathbf{1}}^{(\mathit{mom})}$ inside 
$\wh{\mathcal{A}}_{q}^{(\mathit{mom, ext})}$. 
Therefore, if we describe a bijection between 
the image of $\wh{\mathcal{A}}_{q}^{(\mathit{coor})}$ and 
the image of $\wh{\mathcal{A}}_{\mathbf{1}}^{(\mathit{coor})}$ on one side, 
and a bijection between the image of 
$\wh{\mathcal{A}}_{q}^{(\mathit{mom})}$ and 
$\wh{\mathcal{A}}_{\mathbf{1}}^{(\mathit{mom})}$, 
then a map $\wh{\mathcal{A}}_{\mathbf{1}}^{(\mathit{coor})} \to 
\wh{\mathcal{A}}_{\mathbf{1}}^{(\mathit{mom})}$ induces a map 
$\wh{\mathcal{A}}_{q}^{(\mathit{coor})} \to 
\wh{\mathcal{A}}_{q}^{(\mathit{mom})}$. 

To formalize this, let us first look at the algebra which is obtained 
from $\wt{\mathcal{A}}_{q}$ (generators $p_i$, $x_{\alpha}$, and $h_{\alpha, i}$) by 
adjoining the elements 
\begin{equation*} 
b_{N} \Longrightarrow x_{1}^{- N_1} \dots x_{s}^{- N_{s}}, \quad 
a_{M} \Longrightarrow p_{1}^{- M_1} \dots p_{s}^{- M_{s}}, 
\end{equation*}
where $N = (N_1 \dots, N_s)$ and $M = (M_1, \dots, M_s)$ are multi-indices 
varying over $\mathbb{Z}_{\geqslant 0}^{s}$. 
Adjoin to it the inverse units of measurements $\varkappa_{i}$ and $\theta_{\alpha}$, 
and their inverses $\varkappa_{i}^{-1}$, $\theta_{\alpha}^{-1}$, $i, \alpha \in [s]$. 
Set 
\begin{equation*} 
\wt{p}_{i} := \varkappa_{i} p_{i}, \quad 
\wt{x}_{\alpha} := \theta_{\alpha} x_{\alpha}, 
\end{equation*}
$i, \alpha \in [s]$. 
Recall, that we write $x^{N} = x_{s}^{N_s} \dots x_{1}^{N_1}$, and 
$p^M = p_{s}^{M_s} \dots p_{1}^{M_1}$ for multi-indices $M, N \in \mathbb{Z}^s$. 
For the quantity $S (x) = \sum_{N} (x^N/ N!) b_{N}$ (where $N$ varies over arbitrary 
finite subset of multi-indices from $\mathbb{Z}_{\geqslant 0}^{s}$), we have 
$S (x) \Longrightarrow x^{\bar 0} p^{\bar 0}$. 
If one wants to perceive it as 
$S (x)= \sum_{N} (\wt{x}^N/ N!) \wt{b}_{N}$, then from 
$\wt{x}_{\alpha} \Longrightarrow x^{\bar 0} p^{\bar 0}$, and from  
\begin{equation*} 
\wt{x}_{\lambda_r} \dots \wt{x}_{\lambda_1} = 
\theta_{\lambda_1} \dots \theta_{\lambda_r} 
x_{\lambda_r} \dots x_{\lambda_1} = 
x_{\lambda_r} \dots x_{\lambda_1} \theta_{\lambda_1} \dots \theta_{\lambda_r}, 
\end{equation*} 
where $r \in \mathbb{Z}_{> 0}$, $\lambda_1, \dots, \lambda_{r} \in [s]$, 
one can see, that it is necessary to put 
\begin{equation*} 
\wt{b}_{N} := \theta_{s}^{- N_{s}} \dots \theta_{1}^{- N_{1}} b_{N}, 
\end{equation*}
for every $N \in \mathbb{Z}_{\geqslant 0}^{s}$. 
Having in mind 
$\wt{S} (p) = \sum_{M} (p^M/ M!) a_{M} = \sum_{M} (\wt{p}^M/ M!) \wt{a}_{M}$, set by analogy 
\begin{equation*} 
\wt{a}_{M} := \varkappa_{s}^{- M_{s}} \dots \varkappa_{1}^{- M_{1}} a_{M}, 
\end{equation*}
where $M \in \mathbb{Z}_{\geqslant 0}^{s}$. 
Finally, since for the natural extension of the bracket on $\wt{\mathcal{A}}_{q}$,  
we have $\langle \wt{p}_i, \wt{x}_{\alpha} \rangle = 
\langle \varkappa_{i} p_{i}, x_{\alpha} \theta_{\alpha} \rangle = 
\varkappa_{i} h_{\alpha, i} \theta_{\alpha}$, set 
\begin{equation*} 
\wt{h}_{\alpha, i} := \varkappa_{i} h_{\alpha, i} \theta_{\alpha},
\end{equation*}
so that $\langle \wt{p}_{i}, \wt{x}_{\alpha} \rangle = \wt{h}_{\alpha, i}$, $i, \alpha \in [s]$. 
Since 
\begin{equation*} 
\wt{p}_i, \wt{x}_{\alpha}, \wt{h}_{\alpha, i}, \wt{a}_{M}, \wt{b}_{N} 
\Longrightarrow x^{\bar 0} p^{\bar 0}, 
\end{equation*}
we are now in the situation of the previous section where everything commutes ($q = \mathbf{1}$). 
Replacing the symbols with there analogs carrying the tildes, the result obtained there 
can be formulated as follows. 
To every quantity of the shape 
\begin{equation*} 
\big\langle \wt{x}_{\lambda_{r + 2}}, \dots \big\langle \wt{x}_{\lambda_2}, \big\langle \wt{x}_{\lambda_1}, 
\sum_{M} \frac{\wt{p}^M}{M!} \wt{a}_{M}
\big\rangle \big\rangle \dots \big\rangle \big|_{\wt{p} \to \bar 0}, 
\end{equation*} 
where $r \in \mathbb{Z}_{> 0}$, 
one has associated a polynomial in 
\begin{equation*} 
\wt{h}_{\alpha, i}, \quad 
\wt{b}_{1_{\alpha} + 1_{\beta}}, \quad 
\wt{a}_{1_{i} + 1_{j}}, \quad 
\big\langle \wt{p}_{k_{r' + 2}}, \dots \big\langle \wt{p}_{k_{2}}, \big\langle \wt{p}_{k_1}, 
\sum_{N} \frac{\wt{x}^N}{N!} \wt{b}_{N} 
\big\rangle \big\rangle \big\rangle \big|_{\wt{x} \to \bar 0}, 
\end{equation*}
where $i, j, \alpha, \beta, k_1, k_2, \dots, k_r \in [s]$, $r' \in \mathbb{Z}_{>0}$.  
If we now go back to the variables without tildes $p_i$, $x_{\alpha}$, $h_{\alpha, i}$, 
$a_{M}$, $b_{N}$, via extracting all the factors $\varkappa_i$, $\theta_{\alpha}$ 
by bringing them, say, in front of the sums using the braiding relations, 
then it will turn out that these factors can be cancelled out, thus yielding 
an explicit formula defining $L_{q}$ for $q$ generic. 

Let us describe the $q$-modification of the formulae. 
Define 
$B_{k_1, k_2, \dots, k_{r + 2}}$, $r \in \mathbb{Z}_{> 0}$, and 
$A_{\lambda_{1}, \lambda_{2}, \dots, \lambda_{r' + 2}}$, $r' \in \mathbb{Z}_{> 0}$, from 
\begin{equation*} 
\begin{gathered} 
\big\langle \wt{p}_{k_{r + 2}}, \dots \big\langle \wt{p}_{k_2}, \big\langle \wt{p}_{k_1}, 
\sum_{N} \frac{\wt{x}^N}{N!} \wt{b}_{N} \big\rangle \big\rangle \dots \big\rangle\big|_{\wt{x} \to \bar 0} = 
\varkappa_{k_1} \varkappa_{k_2} \dots \varkappa_{k_{r + 2}} B_{k_{1}, k_{2}, \dots, k_{r + 2}}, \\
\big\langle
\wt{x}_{\lambda_{r' + 2}}, \dots 
\big\langle 
\wt{x}_{\lambda_{2}}, 
\big\langle 
\wt{x}_{\lambda_{1}}, 
\sum_{M} \frac{\wt{p}^M}{M!} \wt{a}_{M} 
\big\rangle 
\big\rangle \dots \big\rangle \big|_{\wt{p} \to \bar 0} = 
\theta_{\lambda_1} \theta_{\lambda_{2}} \dots \theta_{\lambda_{r' + 2}}
A_{\lambda_1, \lambda_2, \dots, \lambda_{r' + 2}}. 
\end{gathered}
\end{equation*}
Adjoin the symbols 
$u_{N} (\alpha, \beta) \Longrightarrow x^{-N + 1_{\alpha} + 1_{\beta}}$, 
$N \in \mathbb{Z}_{\geqslant 0}^{s}$, $\alpha, \beta \in [s]$, to our 
algebra with the generators 
$p_{i}$, $x_{\alpha}$, $h_{\alpha, i}$, $a_{M}$, $b_{N}$ 
($i$ and $\alpha$ vary over $[s]$, and $M$ and $N$ vary over $\mathbb{Z}_{\geqslant 0}^{s}$). 
Then  
\begin{multline*} 
\big\langle 
\wt{p}_{k_r}, \dots, 
\big\langle \wt{p}_{k_2}, 
\big\langle \wt{p}_{k_1}, 
\sum_{\alpha, \beta \in [s]}
\langle \wt{p}_j, \wt{x}_{\alpha} \rangle \wt{b}_{1_{\alpha} + 1_{\beta}} 
\sum_{N} \frac{\wt{x}^N}{N!} \wt{u}_N (\beta, \omega) 
\big\rangle 
\big\rangle \dots 
\big\rangle \big|_{\wt{x} \to \bar 0} 
= \\ =  
\varkappa_{k_1} \varkappa_{k_2} \dots \varkappa_{k_r} \varkappa_{j} 
\wh{U}_{k_1, k_2, \dots, k_r} (j, \omega) \theta_{\omega}, 
\end{multline*} 
where $\wh{U}_{k_1, k_2, \dots, k_r} (j, \omega)$ is some expression 
not containing the inverse units of measurements $\theta_{*}$, $\varkappa_{*}$. 
We can also consider 
$\sum_{\alpha, \beta \in [s]} 
\langle \wt{p}_{j}, \wt{x}_{\alpha} \rangle \wt{b}_{1_{\alpha} + 1_{\beta}} \wt{u}_{\bar 0} (\beta, \omega) = 
\varkappa_{j} \wh{U}_{\emptyset} (j, \omega) \theta_{\omega}$, 
which, on the other hand, according to the previous section, must be 
$\langle \wt{p}_{j}, \wt{x}_{\omega} \rangle = 
\varkappa_{j} \langle p_{j}, x_{\omega} \rangle \theta_{\omega}$, 
so we have just $\wh{U}_{\emptyset} (j, \omega) = 
\langle p_j, x_{\omega} \rangle = h_{\omega, j}$, $j, \omega \in [s]$. 
Now, considering the equations 
(\ref{eq:U0B0_recursive})
linking the coefficients $B_{*}^{(0)}$ and $U_{*}^{(0)}$ 
in the previous section, define recursively: 
\begin{multline}  
\label{eq:barUbarB_recursive}
\bar U_{k_1, k_2, \dots, k_r} (i, \omega) := 
\varkappa_{i}^{-1} 
\varkappa_{k_r}^{-1} \dots \varkappa_{k_1}^{-1} 
\sum_{j, l \in [s]} \sum_{t = 1}^{r} \sum_{\sigma \in \mathit{Sh}_{r} (t)} 
\big[ \varkappa_{i} \varkappa_{j} \varkappa_{k_{\sigma (1)}} \dots \varkappa_{k_{\sigma (t)}}
\times \\ \times 
\bar B_{i, j, k_{\sigma (1)}, \dots, k_{\sigma (t)}} \big] 
\big( 
\varkappa^{- (1_{j} + 1_{l})} \bar a_{j, l}
\big) 
\big\lbrace 
\varkappa_{k_{\sigma (t + 1)}} \dots \varkappa_{k_{\sigma (r)}} \varkappa_{l} 
\bar U_{k_{\sigma (t + 1)}, \dots, k_{\sigma (r)}} (l, \omega) \theta_{\omega} 
\big\rbrace
\theta_{\omega}^{-1}
\end{multline} 
starting from $\bar U_{\emptyset} (i, \omega) := h_{\omega, i}$, 
where $i, \omega, k_1, k_2, \dots, k_r \in [s]$, $r \in \mathbb{Z}_{>0}$, 
and we denote 
$\mathit{Sh}_{r} (t) := \lbrace \sigma \in S_{r} \, | \, \sigma (1) < \dots < \sigma (t) \text{ and } 
\sigma (t + 1) < \dots < \sigma (r) \rbrace$.  
The factors $\theta_{\omega}$ and $\theta_{\omega}^{-1}$ cancel out immediately, 
and if one brings all the coefficients $\varkappa_{*}$ in front of the sums, 
then they cancel out as well, leaving some coefficients inside the sums, which 
are just some products of $q_{i', j'}$, $i', j' \in [2 s]$ stemming from the braidings. 
Therefore, all $\bar U_{k_1, k_2, \dots, k_r} (i, \omega)$ are just 
polynomials in variables $h_{*}$, $\bar a_{*}$, $\bar b_{*}$, and $\bar B_{*}$ 
(we write $*$ instead of blind indices in the subscripts). 

For the next step, in accordance with 
(\ref{eq:UemptysetUk}), (\ref{eq:Vemptyset}), (\ref{eq:Vnu}), 
one needs to consider the ``tilded expressions'' 
\begin{equation*} 
\wt{V}_{\emptyset} (p; \alpha, \omega)|_{\wt{p} \to \bar 0} = 
\sum_{i, j, \alpha, \beta \in [s]}
\langle \wt{x}_{\alpha}, \wt{p}_{i} \rangle \wt{a}_{1_{i} + 1_{j}} 
\Big\lbrace 
\langle \wt{p}_{j}, \wt{x}_{\alpha} \rangle \wt{b}_{1_{\alpha} + 1_{\beta}} 
\sum_{N} \frac{\wt{x}^N}{N!} \wt{u}_{N} (\beta, \omega)
\Big|_{\wt{x} \to \bar 0}  
\Big\rbrace, 
\end{equation*} 
and 
\begin{multline*} 
\wt{V}_{\nu_1, \dots, \nu_r} (\wt{p}; \alpha, \omega)|_{\wt{p} \to \bar 0} = 
\sum_{k_1, \dots, k_r \in [s]} 
\Big( 
\prod_{i' = 1}^{r} 
\sum_{\beta, i \in [s]}
\wt{b}_{1_{\nu_{i'}} + 1_{\beta}} \langle \wt{x}_{\beta}, \wt{p}_{i} \rangle \wt{a}_{1_{i} + 1_{k_{i'}}} 
\Big) 
\sum_{j, l \in [s]}
\langle \wt{x}_{\alpha}, \wt{p}_{j} \rangle 
\times \\ \times
\wt{a}_{1_{j} + 1_{l}} 
\Big\lbrace 
\langle 
\wt{p}_{k_r}, \dots, \langle \wt{p}_{k_2}, \langle \wt{p}_{k_1}, 
\sum_{\sigma, \rho \in [s]}
\langle \wt{p}_{l}, \wt{x}_{\sigma} \rangle \wt{b}_{1_{\sigma} + 1_{\rho}} 
\sum_{N} \frac{\wt{x}^N}{N!} \wt{u}_{N} (\rho, \omega) 
\rangle \rangle \dots \rangle
\Big|_{\wt{x} \to \bar 0}
\Big\rbrace, 
\end{multline*}
for every $\alpha, \omega, \nu_1, \dots, \nu_r \in [s]$, and every $r \in \mathbb{Z}_{>0}$. 
One can see, that 
\begin{equation*} 
\begin{gathered}
\wt{V}_{\emptyset} (p; \alpha, \omega)|_{\wt{p} \to \bar 0} = 
\theta_{\alpha} \wh{V}_{\emptyset} (\alpha, \beta) \theta_{\omega},\\ 
\wt{V}_{\nu_1, \dots, \nu_r} (p; \alpha, \omega)|_{\wt{p} \to \bar 0} = 
\theta_{\nu_1} \dots \theta_{\nu_r} \theta_{\alpha} 
\wh{V}_{\nu_1, \dots, \nu_r} (\alpha, \omega) \theta_{\omega} 
\end{gathered}
\end{equation*}
where $\wh{V}_{\emptyset} (\alpha, \omega)$ and $\wh{V}_{\nu_1, \dots, \nu_r} (\alpha, \omega)$ 
are expressions not containing $\varkappa_{i}$, $\theta_{\alpha}$, $i, \alpha \in [s]$.  
In fact, $\wh{V}_{\emptyset} (\alpha, \omega) = \sum_{i, j \in [s]} \langle x_{\alpha}, p_{i} \rangle 
a_{1_{i} + 1_{j}} \langle p_{j}, x_{\omega} \rangle$. 
The expressions in the curly brackets have already been analysed above. 
Therefore, one defines 
\begin{equation} 
\label{eq:barVbarUemptyset}
\bar V_{\emptyset} (\alpha, \omega) := 
\theta_{\alpha}^{-1} \sum_{\alpha, \beta, i, j \in [s]} 
\big( \theta_{\alpha} \langle x_{\alpha}, p_{i} \rangle \varkappa_{i} \big) 
\big[ 
\varkappa^{- (1_{i} + 1_{j})} \bar a_{i, j}
\big] 
\big\lbrace
\varkappa_{j} \bar U_{\emptyset} (j, \omega) \theta_{\omega}
\big\rbrace \theta_{\omega}^{-1}, 
\end{equation}
and 
\begin{multline}
\label{eq:barVnubarUk} 
\bar V_{\nu_1, \dots, \nu_r} (\alpha, \omega) := 
\theta_{\alpha}^{-1} \theta_{\nu_r}^{-1} \dots \theta_{\nu_1}^{-1} 
\sum_{k_1, \dots, k_r \in [s]} 
\bigg[
\prod_{i' = 1}^{r}
\sum_{\beta, i \in [s]} 
\Big( 
\theta^{- (1_{\nu_{i'}} + 1_{\beta})} \bar b_{\nu_{i'}, \beta}
\Big) 
\times \\ \times 
\Big[ 
\theta_{\beta}^{-1} \langle x_{\beta}, p_{i} \rangle \varkappa_{i}^{-1}
\Big] 
\Big\lbrace 
\varkappa^{- (1_{i} + 1_{k_{i'}})} \bar a_{i, k_{i'}}
\Big\rbrace
\bigg]
\bigg\lbrace 
\sum_{j, l \in [s]} 
\Big( 
\theta_{\alpha} \langle x_{\alpha}, p_{j} \rangle \varkappa_{j}
\Big)
\times \\ \times 
\Big[ 
\varkappa^{- (1_{j} + 1_{l})} \bar a_{j, l}
\Big] 
\varkappa_{k_1} \dots \varkappa_{k_r} \varkappa_{l} 
\bar U_{k_1, \dots, k_r} (l, \omega) \theta_{\omega} 
\bigg\rbrace
\theta_{\omega}^{-1}, 
\end{multline}
for $\alpha, \omega, \nu_1, \dots, \nu_r \in [s]$, $r \in \mathbb{Z}_{> 0}$. 
After concentrating all $\theta_{*}$ and $\varkappa_{*}$ in front of the sums, 
these generators cancel out 
(leaving a sum of terms having the same ``units of measurement''). 
Evaluating the brackets $\langle p_i, x_{\alpha} \rangle$ and 
$\langle x_{\alpha}, p_{i} \rangle$, $i, \alpha, \in [s]$, 
and expressing the coefficients $\bar U_{\emptyset} (j, \omega)$ and 
$\bar U_{k_1, \dots, k_r} (l, \omega)$ as described in 
(\ref{eq:barUbarB_recursive}), 
we can claim that $\bar V_{\emptyset} (\alpha, \omega)$ and 
$\bar V_{\nu_1, \dots, \nu_{r}} (\alpha, \omega)$ are just some 
polynomials in $h_{*}$, $\bar a_{*}$, $\bar b_{*}$, and $\bar B_{*}$ 
with coefficients given by products of the braiding coefficients $q_{i', j'}$, $i', j' \in [2 s]$.  

To make the third step, 
in analogy with 
(\ref{eq:AzeroVzero}), 
we need to consider the tilded expression 
\begin{multline*} 
\big\langle
\wt{x}_{\lambda_{r}}, \dots 
\big\langle \wt{x}_{\lambda_{1}}, 
\big\langle 
\wt{x}_{\mu_{0}}, 
\big\langle 
\wt{x}_{\lambda_{0}}, 
\sum_{M} \frac{\wt{p}^M}{M!} \wt{a}_{M} 
\big\rangle 
\big\rangle 
\big\rangle
\dots \big\rangle \big|_{\wt{p} \to \bar 0} 
= \\ = 
(- 1)^{r - 1} 
\sum_{\substack{f : \lbrace 0, \dots, r \rbrace \to 
\lbrace 0, \dots, r - 1 \rbrace, \\
0 \leqslant f (j') \leqslant j' - 1, \\
j' = 1, \dots, r 
}}
\sum_{\mu_1, \dots, \mu_r \in [s]} 
\prod_{i' = 0}^{r} 
\wt{V}_{\mu_{f^{- 1} (i')}} (\wt{p}; \lambda_{i'}, \mu_{i'})|_{\wt{p} \to \bar 0}, 
\end{multline*}
where $\lambda_{0}, \mu_{0}, \lambda_{1}, \dots, \lambda_{r} \in [s]$, $r \in \mathbb{Z}_{> 0}$. 
The notation $\wt{V}_{\mu_{f^{- 1} (i')}} (\wt{p}; \lambda_{i'}, \mu_{i'})$ 
in the right-hand side is as follows. 
Let the set $f^{-1} (i')$ (where $i' = 1, 2, \dots, r$) be written as 
\begin{equation*} 
f^{-1} (i') = \big\lbrace (f^{-1} (i'))_1 < (f^{-1} (i'))_{2} < \dots < 
(f^{-1} (i'))_{|f^{-1} (i')|} \big\rbrace, 
\end{equation*} 
where $| \cdot |$ denotes the cardinality of a set. 
We put 
\begin{equation*} 
\wt{V}_{\mu_{f^{-1} (i')}} (\wt{p}; \lambda_{i'}, \mu_{i'}) := 
\wt{V}_{\mu_{(f^{-1} (i'))_1}, \dots, \mu_{(f^{-1} (i'))_{|f^{-1} (i')|}}} 
(\wt{p}; \lambda_{i'}, \mu_{i'}). 
\end{equation*}
Therefore, we define: 
\begin{multline}
\label{eq:Ahash} 
A_{\lambda_{0}, \mu_{0}, \lambda_1, \dots, \lambda_r}^{\#} := 
\theta_{\lambda_r}^{-1} \dots \theta_{\lambda_1}^{-1} 
\theta_{\mu_0}^{-1} \theta_{\lambda_0}^{-1} 
\bigg\lbrace 
(- 1)^{r - 1} 
\sum_{\substack{f : \lbrace 0, \dots, r \rbrace \to 
\lbrace 0, \dots, r - 1 \rbrace, \\
0 \leqslant f (j') \leqslant j' - 1, \\
j' = 1, \dots, r 
}}
\sum_{\mu_1, \dots, \mu_r \in [s]} 
\times \\ \times
\prod_{i' = 0}^{r} 
\Big[ 
\theta_{\mu_{(f^{-1} (i'))_{1}}} \dots 
\theta_{\mu_{(f^{-1} (i'))_{|f^{-1} (i')|}}}
\theta_{\lambda_{i'}} 
\bar V_{   \mu_{(f^{-1} (i'))_{1}}, \dots, (f^{-1} (i'))_{|f^{-1} (i')|} } (\lambda_{i'}, \mu_{i'}) 
\theta_{\mu_{i'}}
\Big] 
\bigg\rbrace, 
\end{multline}
for $\lambda_{0}, \mu_{0}, \lambda_{1}, \dots, \lambda_{r} \in [s]$, $r \in \mathbb{Z}_{> 0}$. 
After bringing all $\theta_{*}$ to the left in front of the summation, one observes, that 
they cancel out leaving in the right-hand side a polynomial with respect to 
$\bar V_{*}$. Its coefficients are just some products of the 
braiding coefficients forming the $2 s  \times 2 s$ matrix $q$.   
Since $\bar V_{*}$ are some polynomials in $h_{*}$, $\bar a_{*}$, $\bar b_{*}$, and $\bar B_{*}$, 
these quantities can also be perceived as polynomials in $h_{*}$, $\bar a_{*}$, $\bar b_{*}$, and $\bar B_{*}$. 
The coefficients are some polynomials in the entries of the matrix $q$. 
In remains tho mention, that in the algebra $\mathcal{T}_{q}^{(\mathit{coor})}$, 
which is generated by $h_{*}$, $\bar a_{*}$, $\bar b_{*}$, and $\bar B_{*}$, 
the quantities 
$A_{*}^{\#}$
satisfy just the same commutation relations (i.e. have the same braidings), 
as the quantities $\bar A_{*}$ in the algebra $\mathcal{T}_{q}^{(\mathit{mom})}$, 
which is generated by $h_{*}$, $\bar a_{*}$, $\bar b_{*}$, and $\bar A_{*}$.  
Therefore, define $L_{q} : \mathcal{T}_{q}^{(\mathit{mom})} \to 
\mathcal{T}_{q}^{(\mathit{coor})}$ on the generators as 
\begin{equation} 
\label{eq:Lq_generators}
h_{\alpha, i} \mapsto h_{\alpha, i}, \quad 
\bar a_{i, j} \mapsto \bar a_{i, j}, \quad 
\bar b_{\alpha, \beta} \mapsto \bar b_{\alpha, \beta}, \quad  
\bar A_{\alpha, \beta, \lambda_1, \dots, \lambda_r} \mapsto 
A_{\alpha, \beta, \lambda_1, \dots, \lambda_r}^{\#}, \quad 
\end{equation} 
for every $\alpha, \beta, i, j, \lambda_1, \dots, \lambda_r \in [s]$, $r \in \mathbb{Z}_{> 0}$, 
and extend it uniquely to an algebra isomorphism. 
\qed

To summarize, 
one starts with an affine quantum space $\mathcal{A}_{q}$
with generators $p_{i}$, $x_{\alpha}$, $i, \alpha \in [s]$. 
After that it is necessary to consider symbols with the braiding behaviour 
\begin{equation*} 
\begin{gathered}
h_{\alpha, i} \Longrightarrow x^{1_{\alpha}} p^{1_{i}}, \quad 
\bar a_{i, j} \Longrightarrow x^{- (1_{i} + 1_{j})}, \quad 
\bar b_{\alpha, \beta} \Longrightarrow p^{- (1_{\alpha} + 1_{\beta})}, \\
\bar A_{\alpha, \beta, \lambda_1, \dots, \lambda_r} \Longrightarrow 
x^{1_{\alpha} + 1_{\beta} + 1_{\lambda_1} + \dots + 1_{\lambda_r}}, \quad 
\bar B_{i, j, k_1, \dots, k_r} \Longrightarrow 
p^{1_{i} + 1_{j} + 1_{k_1} + \dots + 1_{k_r}}. 
\end{gathered}
\end{equation*}
Imposing the relations mentioned above, 
one defines an algebra $\mathcal{T}_{q}^{(\mathit{coor})}$ generated by 
$h_{*}$, $\bar a_{*}$, $\bar b_{*}$, and $\bar B_{*}$, and 
and algebra $\mathcal{T}_{q}^{\mathit{(mom)}}$ generated by 
$h_{*}$, $\bar a_{*}$, $\bar b_{*}$, and $\bar A_{*}$. 

\begin{defn} 
The algebra isomorphism $L_{q}: \mathcal{T}_{q}^{(\mathit{coor})} \overset{\sim}{\to} 
\mathcal{T}_{q}^{(\mathit{mom})}$ 
constructed in the proof of the theorem (formulae 
(\ref{eq:barUbarB_recursive}),  
(\ref{eq:barVbarUemptyset}), 
(\ref{eq:barVnubarUk}), 
(\ref{eq:Ahash}), 
(\ref{eq:Lq_generators})), 
is called a \emph{$q$-Legendre transformation (in a point)}. 
\end{defn}

The formulae defining the inverse isomorphism are constructed in a totally similar 
way and define the inverse $q$-Legendre transformation (in a point).

\section{Discussion.}
Let us now go back to the analogy between 
the semiclassical quantum mechanics and the 
quasithermodynamic 
statistical physics  
discussed in the introduction. 
This paper is a natural continuation of 
\cite{RuugeVanOystaeyen}. 
It is important to stress, 
that in order to construct 
a reasonable $q$-analogue of a classical theory (mechanics or thermodynamics), 
it is \emph{not} enough just to apply the $q$-analysis 
(replacing derivatives with $q$-derivatives, factorials with $q$-factorials, etc.).  
It is necessary to introduce the 
Planck-Boltzmann constants $\hbar \to 0$ or $k_{B} \to 0$ into the theory first. 
The paper \cite{RuugeVanOystaeyen} is focused on $\hbar \to 0$, and starts 
with an investigation of a $q$-analogue of the Weyl quantization map in quantum mechanics. 
In the non-$q$-deformed case this is just a symmetrization map linking 
the classical coordinates $x$ and momenta $p$, with the quantized 
coordinates and momenta $\wh{x}$ and $\wh{p}$. 
Trying to construct a reasonable $q$-analogue of such symmetrization map 
in case $q$ is a $2 s \times 2 s$ matrix of formal variables 
(\ref{eq:qq_assumptions}) 
defining the 
braidings 
(\ref{eq:xixi_rels}) 
on the phase space 
(where $s$ is the number of degrees of freedom), one realizes 
that the Planck constant $\hbar$ should \emph{acquire indices}, 
$\hbar \to \hbar_{i, j}$, $i, j \in [2 s]$, and should have same 
commutation properties as the $q$-commutator 
\begin{equation*}
[\wh{z}_{i}, \wh{z}_{j}]_{q} := \wh{z}_{i} \wh{z}_{j} - q_{j, i} \wh{z}_{j} \wh{z}_{i},  
\end{equation*}
where $\wh{z}_{i} = \wh{p}_{i}$, for $i = 1, 2, \dots, s$, and 
$\wh{z}_{i} = \wh{x}_{i - s}$, for $i = s + 1, s + 2, \dots, 2 s$. 
This is a part of a more general construction which we term 
the \emph{bracketing algebra} (or, also the \emph{epoch\'e algebra}). 
One can introduce ``higher order'' Planck constants 
\begin{equation*}
\hbar_{(i, j), k},\quad 
\hbar_{i, (j, k)},\quad 
\hbar_{((i, j), k), l}, \quad 
\hbar_{(i, (j, k)), l}, \quad 
\hbar_{i, ((j, k), l)}, \quad 
\hbar_{i, (j, (k, l))}, \quad 
\hbar_{(i, j), (k, l)}, \quad \text{etc.} 
\end{equation*}
where $i, j, k, l \in [2 s]$, 
which behave like the corresponding $q$-commutators 
constructed from the ``Planck constants'' of the lower orders.  
More precisely, one obtains an algebra with an \emph{infinite} number 
of generators $\hbar_{\Gamma}$ indexed by leaf-labelled planar binary trees $\Gamma$, 
satisfying the relations 
\begin{equation*} 
\hbar_{\Gamma} \hbar_{\Gamma'} - q_{\Gamma', \Gamma} \hbar_{\Gamma'} \hbar_{\Gamma} = 
\hbar_{\Gamma' \vee \Gamma}, 
\end{equation*}  
for any $\Gamma$, $\Gamma'$, where $q_{\Gamma', \Gamma}$ is a certain naturally defined product of $q_{i, j}$, 
and $\Gamma' \vee \Gamma$ denotes the concatenation of trees ($\Gamma'$ becomes the left branch, $\Gamma$ 
becomes the right branch). 
The set of labels of the leaves is just the symbols 
$\lbrace \bar p_1, \dots, \bar p_s, \bar x_1, \dots, \bar x_s \rbrace$ 
corresponding to the classical coordinates and momenta.  
There are also some conditions generalizing $\hbar_{i, j} = - q_{j, i}^{-1} \hbar_{j, i}$, 
related to the symmetry of a tree $\Gamma$, for details refer to \cite{RuugeVanOystaeyen}. 
The generators $\hbar_{\Gamma}$ corresponding to the labelled trees with just one leaf are 
identified with the ``quantized'' coordinates and momenta. 
Therefore, the higher order generators $\hbar_{\Gamma}$ are just as good 
for the role of dynamical quantities, as $x_{*}$ and $p_{*}$. 
Intuitively, the quantization of the Planck constant $\hbar$ consists in replacing  
\begin{equation} 
\label{eq:hbar_epoche}
\hbar \to \lbrace \hbar_{\Gamma} \rbrace_{\Gamma} = 
\lbrace \hbar_{i}, \hbar_{i, j}, \hbar_{(i, j), k}, \hbar_{(i, j), (k, l)}, \dots \rbrace
\end{equation} 

It is quite remarkable to observe, that even in the limit 
$q \to \mathbf{1}$, where $\mathbf{1}$ is a $2 s \times 2 s$ 
matrix with all entries equal to 1, one still has 
an infinite collection of generators. 
This picture is similar to what sometimes happens in quantum statistical physics.  
Instead of considering the standard creation and annihilation operators 
$\psi^{\pm} (x)$, where $x$ varies over the 1-particle configuration space, 
one may wish to split the particles in pairs, triples, and so forth, 
and to consider formally this subsets as new particles. 
This is sometimes termed (due to V.~P.~Maslov \cite{Maslov2}) 
the \emph{``ultrasecond'' quantization}, and it can be useful  
even for a system containing only one sort of identical particles 
(e.g., one can think of the Cooper pairs of electrons 
in the theory of low-temperature superconductivity).    
In other words, one works in terms of different kinds of creation-annihilations operators 
$\psi^{\pm} (x)$, $\psi^{\pm} (x, x')$, $\psi^{\pm} (x, x', x'')$, etc., 
\begin{equation} 
\label{eq:ultrasecond}
\psi^{\pm} (x) \to 
\lbrace \psi^{\pm} (x), \psi^{\pm} (x, x'), \psi^{\pm} (x, x', x''), \dots \rbrace, 
\end{equation}
where $x$ is a point of the 1-particle configuration space, 
$(x, x')$ is a point of the 2-particle configuration space, 
$(x, x', x'')$ is a point of the 3-particle configuration space, etc. 
At the same time, the generators $\hbar_{\Gamma}$ are indexed not 
by finite sequences of elements in $\lbrace \bar p_1, \dots, \bar p_s, 
\bar x_1, \dots, \bar x_s \rbrace$, but by the 
canonical basis of the \emph{free} Lie algebra generated by  
$\lbrace \bar p_1, \dots, \bar p_s, 
\bar x_1, \dots, \bar x_s \rbrace$. 
The commutation relations are similar to the canonical commutation 
relations in the Weyl algebra, but are not completely the same. 
One obtains an infinite chain of relations 
involving higher and higher orders of generators $\hbar_{\Gamma}$ 
(i.e. bigger and bigger trees $\Gamma$). 
It can be perceived as a collection of relations, 
describing a deformation of a deformation of a deformation ... 
(infinite number of times) of the canonical commutation relations.    
It is quite interesting to mention in this connection the work of M.~Kapranov \cite{Kapranov}, 
where he defines ``fat'' non-commutative manifolds and considers a filtration 
on a non-commutative algebra given by the commutators in order to define a 
kind of ``non-commutative neighbourhood'' of an algebraic variety.  

It is worth to point out, that instead of the ``ultrasecond'' quantization 
it is possible to consider just the \emph{``ultra'' quantization}, replacing the 
quantum mechanical creation-annihilation operators $a_{i}^{\pm}$, 
$i \in [s] = \lbrace 1, 2, \dots, s \rbrace$ ($s$ is the number of classical degrees of freedom), with 
\begin{equation} 
\label{eq:ultra}
a_{i}^{\pm} \to \lbrace a_{i}^{\pm}, a_{i, i'}^{\pm}, a_{i, i', i''}^{\pm}, \dots \rbrace, 
\end{equation}
where $i \in [s]$, $(i, i') \in [s] \times [s]$, $(i, i', i'') \in [s] \times [s] \times [s]$, etc. 
In this context, one 
can perceive the ``ultrasecond'' quantization
as the second ``ultra'' quantization.

If one looks at the analogy between the semiclassical ($\hbar \to 0$) 
wave functions 
\cite{Littlejohn}, 
and the quasitermodynamical 
($k_{B} \to 0$) 
partition functions  
\cite{Rajeev1, Rajeev2}, 
it is natural to expect 
that it is of interest to replace the Boltzmann constant $k_{B}$ 
with an infinite collection of generators of the epoch\'e algebra,  
\begin{equation} 
\label{eq:kB_epoche}
k_{B} \to \lbrace (k_{B})_{\Gamma} \rbrace_{\Gamma} = 
\lbrace 
(k_{B})_{i}, 
(k_{B})_{i, j}, 
(k_{B})_{(i, j), k}, 
(k_{B})_{(i, j), (k, l)}, \dots
\rbrace, 
\end{equation} 
where $\Gamma$ varies over the set of all finite leaf-labelled planar binary trees, 
and the labelling set is just the set of symbols denoting the thermodynamic quantities. 
For the two-dimensional thermodynamic system, this set is 
$\lbrace \bar T, \bar S, \bar p, \bar V \rbrace$ 
(temperature, entropy, pressure, volume). 
The generators $(k_{B})_{i}$ corresponding to the trees with only one leaf are associated to  
the thermodynamic quantities themselves 
(for the two-dimensional thermodynamic system, this is 
temperature $T$, entropy $S$, pressure $p$, volume $V$). 
Note, that the $q$-deformation of thermodynamics that exists in the 
physics literature, is usually aimed at investigating the so-called 
non-extensive Tsallis entropy 
\cite{LavagnoSwamy, Tsallis}
(the latter seems to be quite useful also in economics). 
What we obtain is a little different, since, for example, the $q$-Legendre 
transformation described in the present paper, in our opinion, 
is most naturally perceived precisely in terms of the truncation of the 
thermodynamic epoch\'e algebra at trees with two leaves  
(we restrict ourselves to the generators $(k_{B})_{i}$ and $(k_{B})_{i, j}$). 

As already mentioned, taking the limit $q \to \mathbf{1}$ for the 
mechanical epoch\'e algebra, does \emph{not} bring one back immediately to 
classical mechanics. 
This fact, actually, modifies one's understanding of what the mechanical classical limit should be. 
One needs to consider a central extension of the epoch\'e algebra, 
\begin{equation*} 
\hbar_{\Gamma} \hbar_{\Gamma'} - q_{\Gamma', \Gamma} \hbar_{\Gamma'} \hbar_{\Gamma} = 
\eta_{\mathit{mech}} \hbar_{\Gamma' \vee \Gamma}, 
\end{equation*}
where $\eta_{\mathit{mech}}$ is a central generator. 
The classical limit consists in specializing $q \to \mathbf{1}$ \emph{and} 
going to the associated graded with respect to the $\eta_{\mathit{mech}}$-adic filtration. 
The usual classical mechanical quantities fall in degree zero, 
and the higher degrees contain the semiclassical ``corrections''. 
The same happens with the thermodynamic epoch\'e algebra. 
Going from the statistical physics of equilibrium states to the phenomenological thermodynamics 
is implemented by considering a central extension of the thermodynamic epoch\'e algebra,  
\begin{equation*} 
(k_{B})_{\Gamma} (k_{B})_{\Gamma'} - q_{\Gamma', \Gamma} (k_{B})_{\Gamma'} (k_{B})_{\Gamma} = 
\eta_{\mathit{thermo}} (k_{B})_{\Gamma' \vee \Gamma}, 
\end{equation*} 
where $\eta_{\mathit{thermo}}$ is a central generator. 
The thermodynamic limit consists in specializing $q \to \mathbf{1}$ \emph{and} 
taking the associated graded with respect to the $\eta_{\mathit{thermo}}$-adic filtration. 
The phenomenological thermodynamics corresponds to the degree zero component, and 
the higher degrees contain the quasithermodynamic ``fluctuations''. 
Note, that implicitly the step of going to the associated graded  
is already present in the ``ultra'' quantization picture   
(\ref{eq:ultrasecond}), (\ref{eq:ultra}).  

Essentially, what is suggested in \cite{RuugeVanOystaeyen} and the present paper, 
is that one should ``blow up'' the Planck constant $\hbar$ and the Boltzmann constant $k_{B}$, 
replacing them with the generators of the epoch\'e algebras 
(\ref{eq:hbar_epoche}) and (\ref{eq:kB_epoche}).   
Hopefully, this blowing up of the Planck-Boltzmann constants into an infinite number of pieces 
does not leave the quantum statistical physics in ruins.

\end{document}